\catcode`\@=11

\font\tenmsb=msbm10
\font\sevenmsb=msbm7 \font\fivemsb=msbm5  \newfam\msbfam
\textfont\msbfam=\tenmsb\scriptfont\msbfam=\sevenmsb
\scriptscriptfont\msbfam=\fivemsb

\def\hexnumber@#1{\ifnum#1<10 \number#1\else \ifnum#1=10 A\else\ifnum#1=11
 B\else\ifnum#1=12 C\else \ifnum#1=13 D\else\ifnum#1=14 E\else\ifnum#1=15
 F\fi\fi\fi\fi\fi\fi\fi}
 \def\msb@{\hexnumber@\msbfam}

\mathchardef\hbar="0\msb@7E

\def\Bbb{\ifmmode\let\next\Bbb@\else
\def\next{\errmessage{Use \string\Bbb\space only in math mode}}\fi\next}
\def\Bbb@#1{{\Bbb@@{#1}}} \def\Bbb@@#1{\fam\msbfam#1}

\catcode`\@=\active

\def\CR{\hbox{{$\cal R$}}}

\def\rbiprod{{\cdot\kern-.33em\triangleright\!\!\!<}}
\def\lbiprod{{>\!\!\!\triangleleft\kern-.33em\cdot\, }}

\def\note#1{{}}
\def\R{{\Bbb R}}
\def\N{{\Bbb N}}
\def\C{{\Bbb C}}
\def\Z{{\Bbb Z}}
\def\H{{\Bbb H}}

\def\isom{{\cong}}
\def\eps{{\epsilon}}
\def\tens{\mathop{\otimes}}

 \def\proof{{\sl Proof}}
\def\endproof{$\bullet$}
\def\vect{{\bf t}}
\def\vecu{{\bf u}}

\def\Diff{{\rm Diff}}

\def\extd{{\rm d}}
\def\Ad{{\rm Ad}}

\def\id{{\rm id}}
\def\<{\langle}
\def\>{\rangle}

\def\und#1{{\underline {#1}}}

\def\o{{}_{\scriptscriptstyle(1)}}
\def\t{{}_{\scriptscriptstyle(2)}}
\def\fo{{}_{\scriptscriptstyle(4)}}

\def\th{{}_{\scriptscriptstyle(3)}}

\def\nquad{{\!\!\!\!\!\!}}

\def\equad{\nquad}

\def\align#1{\begin{eqnarray*}#1\end{eqnarray*}}
\def\cmath#1{\[\begin{array}{c} #1 \end{array}\]}
\def\eqn#1#2{\begin{equation}#2\label{#1}\end{equation}}

\documentstyle[11pt,epsf]{article}
\newtheorem{lemma}{Lemma}[section]
\newtheorem{propos}[lemma]{Proposition}

\newtheorem{corol}[lemma]{Corollary}

\begin{document}\baselineskip 22pt

{\ }\qquad  \hskip 4.3in DAMTP/97-141
\vspace{.2in}

\begin{center} {\LARGE QUANTUM AND BRAIDED DIFFEOMORPHISM GROUPS}
\\ \baselineskip 13pt{\ }
{\ } \\ Shahn Majid\footnote{Royal Society University Research
Fellow and Fellow of Pembroke College, Cambridge}\\ {\ }\\
Department of Applied Mathematics \& Theoretical Physics\\
University of Cambridge, Cambridge CB3 9EW, UK\\
www.damtp.cam.ac.uk/user/majid
\end{center}
\begin{center}
December, 1997
\end{center}

\vspace{10pt}
\begin{quote}\baselineskip 13pt
\noindent{\bf Abstract} We develop a general theory of `quantum'
diffeomorphism groups
based on the universal comeasuring quantum group $M(A)$ associated
to an algebra $A$ and its various quotients. Explicit formulae are
introduced for this construction, as well as dual quasitriangular
and braided $R$-matrix versions. Among the examples, we construct
the $q$-diffeomorphisms of the quantum plane $yx=qxy$, and recover
the quantum matrices $M_q(2)$ as those respecting its braided group
addition law.
\end{quote}
\baselineskip 22pt

\section{Introduction}

In \cite{Swe:hop} one finds a standard construction for a measuring
coalgebra $M(A,B)$ as a universal object for coalgebras `acting'
from algebras $A$ to $B$. The diagonal case $M(A,A)$ is a
bialgebra. The construction is the analogue of the classical
automorphism group except now as the universal or `maximal' object
in the category of bialgebras rather than of groups. Until now,
however, this measuring bialgebra construction has been little
studied in the modern quantum groups literature, although see
\cite{Bat:dif}\cite{Wan:sym}. Probably the main reason for this is
the lack of explicit formulae: being defined as a universal object
it is generally hard to compute.

In the present paper we develop a much more explicit and computable
version of this construction, namely a bialgebra $M(A)$ defined
directly by the structure constants of an associative algebra $A$,
which we call the `comeasuring bialgebra'. It is the
arrows-reversed version of the standard construction but, by
contrast, is defined by generators and relations in a familiar way.
It has many applications and, in fact, a quasi-quantum group
version is used in our forthcoming paper \cite{AlbMa:qua} with H.
Albuquerque as a definition of the automorphism object of
quasi-associative algebras such as the octonions. By contrast, here
we study in much more detail the strictly associative setting and
its `geometrical' applications. We also generalise the construction
to the general braided, where $A$ is an algebra in a braided
category such as that associated to a Yang-Baxter matrix. Such
algebras abound in the theory of $q$-deformations as the
`geometrical objects' on which quantum groups act as symmetries.

The automorphism object $M(A)$ clearly plays the role of
`diffeomorphisms' in a noncommutative geometry setting (where an
algebra $A$ is viewed as like the `functions' on some space). This
is our point of view, and we will use the corresponding terminology
throughout the paper. We also consider briefly quotients of $M(A)$
that preserve a given differential structure on $A$, but for the
most part we work with the universal differential calculus
canonically associated to an algebra. Also, we use the terms
`quantum group' and `braided group' a little loosely, without
requiring the existence of an antipode or `group inversion'. In
examples, the restricted comeasuring bialgebras $M_0(A)$ will often
have an antipode or will admit one by adjoining inverse
determinants etc. Within these limitations, we provide a general
approach to quantum diffeomorphisms which includes
$q$-diffeomorphisms of the line (i.e. some version of a q-Virasoro
quantum group) and of the quantum plane, as well as of
finite-dimensional algebras. Our approach is further justified by
showing that elements of `quantum geometry' may be built around
these objects, along the lines of \cite{BrzMa:gau}\cite{Ma:rie}.
The notion of quantum diffeomorphism group should also be viewed as
a step towards the notion of `quantum manifold', which is a
long-term motivation for the present work.

An outline is the following. In the preliminary Section~2 we
formulate the arrows-reversed measuring bialgebra construction  and
obtain explicit formulae for it and its natural quotients. Apart
from the forthcoming paper \cite{AlbMa:qua} in the quasiassociative
setting, these formulae appear to be new. In Section~3 we compute
several examples of comeasuring algebras, exploring their role as
`quantum diffeomorphisms' of polynomial and discrete spaces.
Section~3.7 contains the maximal comeasuring bialgebra of the
quantum plane, while Section~3.8 obtains the $2\times 2$ quantum
matrices $M_q(2)$ as the q-diffeomorphisms respecting addition.
Note that this result is somewhat different from the
characterisation of quantum matrices in \cite{Man:non}, which is
based on a certain construction of `endomorphisms' for quadratic
algebras and their duals. In Section~4 we turn to general
$R$-matrix constructions. We give a dual quasitriangular version
$M(R,A)$ of the comeasuring construction, including the
applications to the line and the quantum plane. Also using
R-matrices in this section are braided group versions of our
constructions.

\section{General constructions}

The abstract definition of the comeasuring bialgebra is obtained
from \cite{Swe:hop} by reversing arrows giving a universal
comeasuring  algebra $M(A,B)$ `coacting' from algebras $A$ to $B$.
Here the arrows of the coacted-upon objects $A,B$ are not reversed,
i.e. we leave these as algebras and do not make them into
coalgebras as a full dualisation would do. Also, we concentrate on
the diagonal case $M(A)=M(A,A)$ since we will not have much to say
about the general non-diagonal case. Their formulation is, however,
strictly analogous. We work over a general field $k$.

Thus, by definition, a {\em comeasuring} of a unital algebra $A$ is
a pair $(B,\beta)$ where $B$ is a unital algebra and $\beta:A\to
A\tens B$ is an algebra map to the tensor product algebra. We
define $(M(A),\beta_U)$, when it exists, to be the initial object
in the category of  comeasurings of $A$, i.e. a comeasuring such
that for any $(B,\beta)$ there exists a unique algebra map
$\pi:M(A)\to B$ such that $\beta=(\id\tens\pi)\beta_U$.

\begin{propos}cf\cite{Swe:hop} $M(A)$, when it exists, is a
bialgebra and $\beta_U$ is a coaction of it on $A$ as an algebra.
Any other coaction of a bialgebra on $A$ as an algebra is a
quotient of this one.
\end{propos}
\proof This is elementary. We note that $M(A)\tens M(A),
(\beta_U\tens\id)\circ\beta_U$ is also a comeasuring. Hence there
is an algebra map $\Delta:M(A)\to M(A)\tens M(A)$ and
$(\beta_U\tens\id)\circ\beta_U=(\id\tens\Delta)\circ\beta_U$. It
remains to show that $\Delta$ is coassociative. For this, consider
$M(A)^{\tens
3},(\beta_U\tens\id\tens\id)\circ(\beta_U\tens\id)\circ\beta_U$ as
another comeasuring. The map $\pi:M(A)\to M(A)^{\tens 3}$ in this
case is such that $(\id\tens\pi)\circ\beta_U$ is the comeasuring
map associated to $M(A)^{\tens 3}$. Both
$(\Delta\tens\id)\circ\Delta$ and $(\id\tens\Delta)\circ\Delta$
clearly fulfill the role of $\pi$, and since $\pi$ is unique, these
maps coincide. Finally, $k,\beta(a)=a\tens 1$ is a comeasuring and
$\eps:M(A)\to k$ is the induced map. It is easy to see that it
provides a counit. Given any other coaction of a bialgebra $B$ on
$A$ as an algebra (i.e. $A$ a $B$-comodule algebra), the fact that
$B$ comeasures gives the required map $\pi:M(A)\to B$. \endproof

When $A$ is nonunital, we can follow the same definitions while
omitting the conditions that $B$ is unital and that $\beta,\pi$
respect the unit. In this case it is clear that the universal
  object $M(A)$ is a not-necessarily unital
bialgebra. We call it the `nonunital version' of the comeasuring
bialgebra. We can, however, always formally adjoin a unit to it,
extending $\Delta,\eps$ by $\Delta(1)=1\tens 1$ and $\eps(1)=1$. We
denote this extension by $M_1(A)$.

Suppose now  that $A$ is finite-dimensional and let $\{e_i\}$ be a
basis. We let $e_ie_j=c_{ij}{}^k e_k$ define its structure
constants.

\begin{propos} $M_1(A)$ is generated by $1$
and a matrix $\vect=(t^i{}_j)$ of generators, with relations and
coproduct
\[ c_{ij}{}^a t^k{}_a=c_{ab}{}^kt^a{}_it^b{}_j,\quad
\Delta t^i{}_j=t^i{}_a\tens t^a{}_j,\quad \eps(t^i{}_j)=\delta^i{}_j.\]
The map $\beta_U(e_i)=e_a\tens t^a{}_i$ is the coaction.
\end{propos}
\proof It is easy to verify that $\Delta$ extends as an algebra map
and that $\beta_U$ is a coaction. For the former, the proof is
\[ \Delta(c_{ab}{}^kt^a{}_it^b{}_j)=c_{ab}{}^kt^a{}_m t^b{}_n
\tens t^m{}_i t^n{}_j=
c_{mn}{}^at^k{}_a\tens t^m{}_it^n{}_j=c_{ij}{}^b t^k{}_a\tens
t^a{}_b=\Delta(c_{ij}{}^a t^k{}_a).\]
 (This is also a special case
of the quasi-Hopf algebra construction in \cite{AlbMa:qua} or of
the braided case in Section~4). Now, let $(B,\beta)$ be a
comeasuring and define $\pi(t^i{}_j)\in B$ by
$\beta(e_i)=e_a\tens\pi(t^a{}_i)$. Then $\beta$ is an algebra map
is the assertion that $c_{ij}{}^a
\pi(t^k{}_a)=c_{ab}{}^k\pi(t^a{}_i)\pi(t^b{}_j)$, i.e. $\pi:M(A)\to B$
defined in this way extends as an algebra map. Thus $M(A)$ as
generated by the $t^i{}_j$ has the required universal property.
When we adjoin $1$, we obtain $M_1(A)$ as stated.
\endproof

Note that the explicit proof that $\Delta$ provides a bialgebra
works similarly for {\em any} tensor $c_{i_1\cdots
i_m}{}^{j_1\cdots j_n}$ of rank $(m,n)$, i.e. one has an associated
bialgebra with the matrix coalgebra and the relations
\eqn{Mtens}{c_{a_1\cdots a_m}{}^{i_1\cdots i_n} t^{a_1}{}_{j_1}
\cdots t^{a_m}{}_{j_m}
=t^{i_1}{}_{a_1}\cdots t^{i_n}{}_{a_n} c_{j_1\cdots j_m}{}^{a_1\cdots a_n}.}
The associated bialgebra may, however, be trivial. In our case the
associativity of $A$, which is the equation
$c_{ij}{}^ac_{ak}{}^l=c_{ia}{}^l c_{jk}{}^a$, is needed for the
interpretation as universal comeasuring object. It also tends to
ensure that the ideal generated by the stated relations is not too
big. Specifically,
\cmath{c_{ab}{}^c (t^a{}_i t^b{}_j)t^d{}_k c_{cd}{}^l=c_{ij}{}^b
t^c{}_b t^a{}_k c_{ca}{}^l
=c_{ij}{}^b c_{bk}{}^a t^l{}_a
=c_{ib}{}^a t^l{}_a c_{jk}{}^b\\
=c_{ac}{}^l t^a{}_i t^c{}_b c_{jk}{}^b= c_{ac}{}^l c_{bd}{}^c t^a{}_i
(t^b{}_j t^d{}_k)
=c_{ab}{}^c  t^a{}_i (t^b{}_j t^d{}_k)c_{cd}{}^l}
holds automatically in $M_1(A)$.

Now suppose that $A$ is unital and $e_0=1$ is a basis element. We
let $e_i$, $i=1,\cdots,\dim(A)-1$ be the remaining basis elements.

\begin{propos} The comeasuring bialgebra $M(A)$ is the quotient of
 $M_1(A)$ by $t^0{}_0=1,t^i{}_0=0$. Explicitly, it has the
relations
\[ c_{ij}{}^a t^k{}_a=c_{ab}{}^kt^a{}_it^b{}_j+b_it^k{}_j+t^k{}_i b_j,
\quad
c_{ij}{}^0+c_{ij}{}^ab_a=c_{ab}{}^0t^a{}_it^b{}_j+b_ib_j\] and
coproduct, coaction
\cmath{ \Delta t^i{}_j=t^i{}_a\tens t^a{}_j,\quad \Delta b_i=b_a\tens
t^a{}_i+1\tens b_i,
\quad\eps(t^i{}_j)=\delta^i{}_j,\quad\eps(b_i)=0\\
\beta_U(e_i)=1\tens b_i+e_a\tens t^a{}_i,} where $b_i=t^0{}_i$.
\end{propos}
\proof Here $t^0{}_0-1,t^i{}_0$ generates a biideal of $M_1(A)$ and
hence provides a quotient bialgebra. For a direct proof that
$\Delta$ as stated extends to products as an algebra map, we have
\align{&&\equad \Delta c_{ij}{}^a t^k{}_a=c_{ij}{}^at^k{}_b\tens t^b{}_a=
c_{ab}{}^c t^k{}_c\tens t^a{}_it^b{}_j+t^k{}_c\tens
b_it^c{}_j+t^c{}_i b_j\\ &&=c_{cd}{}^kt^c{}_at^d{}_b\tens
t^a{}_it^b{}_j+b_at^k{}_b\tens t^a{}_it^b{}_j+t^k{}_a b_b\tens
t^a{}_it^b{}_j+t^k{}_c\tens b_it^c{}_j+t^c{}_i b_j\\
&&=\Delta(c_{ab}{}^kt^a{}_it^b{}_j+b_it^k{}_j+t^k{}_i b_j)} using
the relations and coproduct stated for $M(A)$. Similarly,
\align{&&\equad\Delta(c_{ab}{}^0t^a{}_it^b{}_j+b_ib_j)=c_{ab}{}^0t^a{}_c
t^b{}_d\tens t^c{}_i t^d{}_j+b_ab_b\tens t^a{}_i t^b{}_j+ 1\tens
b_ib_j+b_a\tens t^a{}_i t^b{}_j+b_a\tens b_i t^a{}_j\\ &&=
 1\tens c_{ab}{}^0 t^a{}_i t^b{}_j+1\tens b_ib_j+b_a\tens c_{ij}{}^b
t^a{}_b= c_{ij}{}^01\tens 1+c_{ij}{}^a1\tens b_a+c_{ij}{}^bb_a\tens
t^a{}_b=\Delta(c_{ij}{}^0+c_{ij}{}^ab_a)} for the second set of
relations. The coproduct has a standard form which is clearly
coassociative, hence we obtain a bialgebra. It inherits the
coaction as shown. Conversely, suppose that $B,\beta$ is a
comeasuring and let $\pi(b_i),\pi(t^i{}_j)\in B$ be defined by
$\beta(e_i)=1\tens
\pi(b_i)+e_a\tens
\pi(t^a{}_i)$.  Similarly to the nonunital case, it follows from $\beta$
a unital algebra map that $\pi$ extends as a unital algebra map and
$\beta=(\id\tens\pi)\circ\beta_U$ by construction. \endproof

Note that these constructions are independent of the choice of
basis (beyond the choice $e_0=1$ in the unital case) because they
are abstractly defined. In a new basis $e_i'=e_a\Lambda^a{}_i$, the
generators of $M_1(A)$ are
\eqn{M1trans}{ t'{}^i{}_j=\Lambda^{-1}{}^i{}_a t^a{}_b \Lambda^b{}_j.}
In the unital case we require that the transformation is of the
form $e'_0=e_0$ and $e'_i=e_0\lambda_i+e_a\Lambda^a{}_i$ where the
indices on $\lambda,\Lambda$ run from $1,\cdots,\dim(A)-1$.  Then
the transformed generators of $M(A)$ are
\eqn{Mtrans}{ b'{}_i=b_a\Lambda^a{}_i+\lambda_i-\lambda_b
\Lambda^{-1}{}^b{}_c t^c{}_d\Lambda^d{}_i,
\quad t'^i{}_j=\Lambda^{-1}{}^i{}_at^a{}_b\Lambda^b{}_j.}

If we are given slightly more structure, namely a linear splitting
$A=1\oplus A'$ then we can define a restricted comeasuring
bialgebra $M_0(A)$ as the universal object for comeasurings that
respect the splitting, i.e. such that $\beta(A')\subset A'\tens B$
in addition to $\beta(1)=1$ as before.

\begin{propos} $M_0(A)$ is the quotient of $M(A)$ by $b_i=0$, i.e. it is
the
associative algebra generated by $1,t^i{}_j$ where $i,j=1,\cdots
\dim(A)-1$, and the relations
\[ c_{ij}{}^a t^k{}_a=c_{ab}{}^kt^a{}_it^b{}_j,\quad
c_{ij}{}^0=c_{ab}{}^0t^a{}_it^b{}_j\] and matrix coalgebra. Its
coaction is $\beta_U(e_i)=e_a\tens t^a{}_i$.
\end{propos}
\proof We suppose that a splitting $A=1\oplus A'$ is given and $\{e_i\}$
are a
basis of $A'$. We clearly have a biideal of $M(A)$ generated by the
$b_i$ and hence a quotient bialgebra. It inherits the relations and
coaction shown. Conversely, if $(B,\beta)$ is a comeasuring
preserving $A'$, we define $\pi(t^i{}_j)\in B$ by
$\beta(e_i)=e_a\tens\pi(t^a{}_i)$ and $\pi(1)=1$. \endproof

This recovers the associative case of the formulae used in
\cite{AlbMa:qua}, but now as a universal comeasuring object. It
depends on the basis only through the splitting of $1$ defined by
it, i.e. is independent of the basis of $A'$ (one may make
transformations as above with $\lambda=0$). In summary, associated
to a split unital algebra $A$, we have a sequence of bialgebras
\[ M_1(A)\to M(A)\to M_0(A)\]
where $M(A)$ is the `usual' definition dual to the construction in
\cite{Swe:hop}, respecting the unit but not its splitting, $M_1(A)$
is the unital extension of the nonunital version, and $M_0(A)$ is
the version respecting the unit and in addition its splitting.

We have similar formulae for $M_1(A,B),M(A,B),M_0(A,B)$ as algebras
of `maps' between algebras. For example, if $B$ has structure
constants $d_{\alpha\beta}{}^\gamma$ in basis $\{f_\alpha\}$, then
$M_1(A,B)$ is generated by $1$ and a rectangular matrix
$\{t^i{}_\alpha\}$ of generators with relations
\[ c_{ij}{}^k t^i{}_\alpha t^j{}_\beta=d_{\alpha\beta}{}^\gamma
t^k{}_\gamma.\]
Similarly for its quotients. Also, although we have assumed finite
bases for our explicit formulae, similar formulae hold more
generally with countable basis. In this case one may have to work
formally or, more precisely, with a suitable completion of tensor
products. In this respect, the original measuring (rather than
comeasuring) construction in \cite{Swe:hop} is better behaved. On
the other hand, we have the advantage in the dual formulation of
explicit formulae as algebras with generators and relations.

\section{Quantum diffeomorphisms of polynomial algebras}

We now compute the comeasuring bialgebras for polynomials $\C[x]$
and their quotients, and justify their role as `diffeomorphisms'.
Further justification is in Section~5 where we consider the
preservation of nonuniversal differential calculi; for the moment
the underlying calculus is the universal one canonically associated
to the algebra. In the case of $\C[x]$, the full comeasuring
bialgebras $M_1(\C[x])$ and $M(\C[x])$ require formal powerseries
and are included for motivation only, but $M_0(\C[x])$ is a
completely algebraic object.

Clearly, the algebraic version of a diffeomorphisms $\C[x]\to
\C[x]$ is a polynomial map $x\mapsto
a(x)=a_0+a_1x+a_2x^2+\cdots$, i.e. a polynomial. If we compose two
such diffeomorphisms $a,b$, then
\[ a\circ b(x)=\sum_j a(x^j)b_j
=\sum_j\sum_{n_1,\cdots,n_j}x^{n_1+\cdots+n_j}a_{n_1}\cdots a_{n_j}b_j
=\sum_i x^i \sum_j \sum_{n_1+\cdots+n_j=i}a_{n_1}\cdots a_{n_j} b_j\]
where $i,j,n_1\in\Z_+$, etc. Here $\Z_+$ denotes the natural
numbers including 0. On the other hand, there is not necessarily a
polynomial inverse. Therefore the `coordinate ring' of the
semigroup of such diffeomorphisms is the commutative bialgebra
$\Diff(\C[x])=\C[t_i|\ i\in\Z_+]$ with countable generators $t_i$
and the coaction
\eqn{coactpoly}{ x\mapsto \sum_i x^i \tens t_i}
as a formal power series. The corresponding coproduct and counit is
\eqn{deltapoly}{ \Delta t_i=\sum_j \sum_{n_1+\cdots+n_j
=i}t_{n_1}\cdots t_{n_j}\tens t_j,
\quad \eps(t_i)=\delta_{i,1}}
and also involves powerseries. On the other hand, if we restrict
attention to diffeomorphisms which preserve the point zero in the
line, we consider only polynomials $a(x)=a_1x+a_2x^2+\cdots$, i.e.
without constant term. This restricted diffeomorphism bialgebra is
the commutative bialgebra $\Diff_0(\C[x])=\C[t_i|\ i\in\N]$ with
the same formulae as above but with $i,j,n_1\in \N$ etc. Here $\N$
denotes the natural numbers not including $0$. In this case
(\ref{deltapoly}) has only a finite number of terms. Thus the
algebraic diffeomorphisms fixing 0 form a bialgebra with algebraic
coalgebra structure and countable generators.

\subsection{Comeasurings of the line}

We now apply the constructions in Section~2, and recover
(noncommutative) versions of the above diffeomorphism bialgebras as
$M(\C[x])$ and $M_0(\C[x])$ respectively. Thus, we take $A=\C[x]$
and basis $e_0=1$ and $e_i=x^i$ for $i\in\N$. Then the structure
constants are
\eqn{conpoly}{ c_{ij}{}^k=\delta_{i+j}^k,\quad c_{ij}{}^0=0}
and the comeasuring bialgebra $M(A)$ from Proposition~2.3 has
generators $1, b_i$ and $t^i{}_j$ with relations
\eqn{fullpoly}{ b_{i+j}=b_ib_j,\quad t^k{}_{i+j}
=\sum_{m+n=k}t^m{}_i t^m{}_j+b_it^k{}_j+t^k{}_i b_j}
for all indices in $\N$. These relations allow us to consider
$t_0\equiv b_1$ and $t_i\equiv t^i{}_1$ as the generators,
obtaining the other $b_i,t^i{}_j$ inductively as
\eqn{tpoly}{ t^i{}_j=\sum_{n_1+\cdots+n_j=i}t_{n_1}\cdots t_{n_j},
\quad b_i=(t_0)^i}
for $i,j\in \N$ and $n_1\cdots n_j\in\Z_+$. In this case the
relations (\ref{fullpoly}) become empty, so $M(\C[x])=\C\<t_i|\
i\in\Z_+\>$, the free algebra on countable generators. Its
coproduct and coaction from Proposition~2.3 therefore takes the
form (\ref{deltapoly}) and (\ref{coactpoly}) with all indices from
$\Z_+$, i.e. $M(\C[x])$ has the same form as $\Diff(\C[x])$ except
that the generators $t_i$ are totally noncommuting.

The restricted comeasuring bialgebra $M_0(\C[x])$ is the quotient
of this where we set $t_0=0$. Equivalently, working from
Proposition~2.4, it has generators $t^i{}_j$ for $i,j\in\N$ with
relations
\eqn{relpoly}{ \sum_{m+n=k} t^m{}_i t^n{}_j=t^k{}_{i+j}}
where all indices are in $\N$.  As before, this implies
\eqn{tnpoly}{ t^i{}_j=\sum_{n_1+\cdots+n_j=i}t_{n_1}\cdots t_{n_j}}
where, $t_i\equiv t^i{}_1$ and now the $n_1,\cdots,n_j\in\N$. This
implies in particular that $t^i{}_j=0$ for $j>i$, which ensures
that the matrix coproduct is now a finite sum of terms. In this
way, $M_0(\C[x])=\C<t_i|\ i\in\N>$, the free algebra, with
coproduct (\ref{deltapoly}) and coaction (\ref{coactpoly}) where
all indices are in $\N$. Thus, $M_0(\C[x])$ has the same form as
$\Diff_0(\C[x])$ except that its generators are totally
noncommuting.

Finally, the extended comeasuring bialgebra $M_1(\C[x])$ has no
classical analogue (since usual diffeomorphisms preserve the
constant function $1\in\C[x]$). From Proposition~2.2 it is
generated by $1$ and $t^i{}_j$, where $i,j\in\Z_+$ and relations
(\ref{relpoly}) (but now with all indices in $\Z_+$). This time we
can reduce the matrix generators to two sequences $E_i\equiv
t^i{}_0$ and $t_i\equiv t^i{}_1$ for $i\in\Z_+$. The others are
recovered by (\ref{tpoly}) where now $i,n_1,\in\Z_+$ etc., and
$j\in\N$. In this way we find that $M_1(\C[x])$ is generated by
$1,E_i,t_i$ for $i\in\Z_+$, with the residual relations and
coalgebra
\[ \sum_{m+n=i}E_m E_n=E_i,\quad \sum_{m+n=i}E_m t_n=t_i=\sum_{m+n=i}t_mE_n\]
\[ \Delta E_j=\sum_i E_i\tens \sum_{n_1+\cdots +n_j=i}t_{n_1}\cdots t_{n_j},\quad
\Delta t_i=\sum_j \sum_{n_1+\cdots +n_j=i}t_{n_1}\cdots t_{n_j}\tens t_j\]
and $\eps(E_i)=\delta^i{}_0$, $\eps(t_i)=\delta^i{}_1$, where all
indices are from $\Z_+$. We note that the constrained sums in these
expressions are the usual convolution product $*$ on sequences, and
in that notation the relations and coproduct of $M_1(\C[x])$ are
\[ E*E=E,\quad E*t=t=t*E,\quad \Delta E_j=\sum_i E_i\tens (t*\cdots
*t)_i,\quad \Delta t_i=\sum_j (t*\cdots*t)_i\tens t_j\]
where the $*$-products are $j$-fold. Finally, the coaction is
\[ 1\mapsto\sum_i x^i\tens E_i,\quad x\mapsto \sum_i x^i\tens t_i.\]
The $M(\C[x])$ above is quotients of this non-unit-preserving
diffeomorphism group by setting $E_i=\delta^i{}_0$.

\subsection{Comeasurings of Grassmann and anyonic variables}

Here we consider `diffeomorphisms' of the finite-dimensional
quotients $A=\C[x]/x^N=0$. The case $N=2$ is where $x$ is a
fermionic or `Grassmann' variable. The case $N>2$ is an anyonic
variable. The addition law and geometry of this `anyonic line' can
be found in \cite{Ma:book}.

We take basis $e_i=x^i$ where $i=1,\cdots,N-1$ and $e_0=1$. Then
$c_{ij}{}^k$ have the same form (\ref{conpoly}) but restricted in
range to $0<i,j,k<N$. In this case  $M=\C<t_i|\ i=0,\cdots,N-1>$ is
the free algebra and has coproduct (\ref{deltapoly}) with indices
including $0$. The restricted comeasuring bialgebra is the free
algebra $M_0=\C<t_i|\ i=1,\cdots,N-1>$ with coproduct likewise from
(\ref{deltapoly}) but with all indices now excluding $0$. On the
other hand the comeasuring bialgebra is now finitely generated and
does not need any formal powerseries. (One can then obtain
$M(\C[x])$ in the projective limit $N\to\infty$.)

For the fermionic case $x^2=0$, the comeasuring bialgebra is
$M=\C<b,t>$ with
\[ \Delta b=1\tens b+b\tens t,\quad \Delta t=t\tens t,\quad \eps(b)=0,
\quad \eps(t)=1,\]
i.e. the matrix of generators has the form $\pmatrix{1&b\cr0&t}$.
The restricted comeasuring bialgebra is just $M_0=\C[t]$ with
$\Delta t=t\tens t$ and $\eps(t)=1$. The extended comeasuring
bialgebra $M_1$ is generated by $1$ and $\pmatrix{a&b\cr c&d}$ with
the relations
\[a^2=a,\quad ac+ca=c,\quad ab=b=ba,\quad ad+cb=d=da+bc\]
and the usual matrix coalgebra.

For an anyonic variable with $x^3=0$, the comeasuring bialgebra is
$M=\C<b,t,s>$ with
\cmath{\Delta b=1\tens b+b\tens t+b^2\tens s,\quad \Delta t=t\tens
t+(bt+tb)\tens s\\
\Delta s=s\tens t+(t^2+sb+bs)\tens s,\quad \eps(b)=\eps(s)=0,\quad
\eps(t)=1,}
i.e. the matrix of generators takes the form $\pmatrix{1&b&b^2\cr
0&t& bt+tb \cr 0&s& t^2+sb+bs}$. The restricted comeasuring
bialgebra is $M_0=\C<t,s>$ with
\[ \Delta t=t\tens t,\quad \Delta s=s\tens t+t^2\tens s,\quad\eps(t)=1,
\quad\eps(s)=0,\]
i.e. the matrix of generators take the form $\pmatrix{t&0\cr
s&t^2}$.

\subsection{Comeasurings of roots of unity}

For another class of discrete spaces (which we consider as discrete
models of a circle), we consider the algebra $A=\C
\Z_N=\C[x]/x^N=1$. We take basis $e_0=1$ and $e_i=x^i$ as before, with
\[ c_{ij}{}^k=\delta^k{}_{i+j},\quad c_{ij}{}^0=\delta_{i,-j},\]
where $0<i,j,k<N$ are added modulo $N$. One may put these into
Propositions~2.2-2.4 to obtain the comeasuring bialgebras.

Actually, it is easier in this class of examples to compute $M_1$
first and then quotient it. This has generators $1,t^i{}_j$ with
relations (\ref{relpoly}), where $i,j\in\Z_N$ including zero. The
convolution formula (\ref{tpoly}) now applies for all $i,j$ since
$\Z_N$ is a group, so generators may all be written in terms of
$t_i\equiv t^i{}_1$ for $i\in\Z_n$. In particular,
\[ E_i\equiv t^i{}_0=\sum_{n_1+\cdots+n_N} t_{n_1}\cdots t_{n_N},\]
where all indices are in $\Z_N$. Thus, the extended comeasuring
bialgebra $M_1$ is generated by $1,t_i$ for $i\in\Z_N$ with the
relations
\[ \sum_{n_1+\cdots +n_{N+1}=i}t_1\cdots t_{N+1}=t_i,\quad\forall
i\in\Z_N\]
and coproduct (\ref{deltapoly}) with all indices in $\Z_N$.

The comeasuring algebra $M$ is given by setting $E_i=\delta^i{}_0$.
Hence it is generated by $1,t_i$ for $i\in\Z_N$ modulo the
relations
\[ \sum_{n_1+\cdots+n_N=i}t_{n_1}\cdots t_{n_N}=\delta^i{}_0,\quad
\forall i\in\Z_N\]
and the same form of coproduct. Here
\[ b_i\equiv t^0{}_i=\sum_{n_1+\cdots+n_i=0}t_{n_1}\cdots t_{n_i},\quad
t^i{}_j=\sum_{n_1+\cdots+n_j=i}t_{n_1}\cdots t_{n_j}\] for
$i,j=1,\cdots,N-1$ and $n_1\in\Z_N$, etc.,   provide the generators
in the form of Proposition~2.3.

Finally, the restricted comeasuring bialgebra $M_0$ is given by
setting $b_i=0$. Hence $M_0$ is generated by $1,t_i$ where
$i=1,\cdots,N-1$ modulo the  relations
\[\sum_{n_1+\cdots+n_N=0}t_{n_1}\cdots t_{n_N}=1,\quad
\sum_{n_1+\cdots+n_N=i}t_{n_1}\cdots t_{n_N}=0,
\quad \sum_{n_1+\cdots n_i=0}t_{n_1}\cdots t_{n_i}=0\]
for $i,n_1$, etc., are in the range $1,\cdots,N-1$ and addition
modulo $\Z_N$. The coproduct is (\ref{deltapoly}) with all indices
similarly in this range.

For $x^2=1$, the extended comeasuring bialgebra $M_1$ is generated
by $1,b,t$ modulo the relations and coalgebra
\[ (b\pm t)^3=(b\pm t),\quad \Delta t=t\tens t+(bt+tb)\tens b,\quad
\Delta b=(b^2+t^2)\tens b+b\tens t,\quad
\eps(t)=1,\quad \eps(b)=0.\]
The matrix of generators has the form $\pmatrix{b^2+t^2&b\cr
bt+tb&t}$. The comeasuring bialgebra $M$ is obtained by setting
$b^2+t^2=1$, $bt+tb=0$ and is therefore generated by $1,b,t$ modulo
the stronger relations and coalgebra
\[ (b\pm t)^2=1,\quad \Delta t=t\tens t,\quad \Delta b=1\tens b
+b\tens t,\quad
\eps(t)=1,\quad \eps(b)=0.\]
Although not a Hopf algebra, there is a formal antipode
\[ St={t\over 1-b^2},\quad Sb=-{bt\over 1-b^2}\]
as a powerseries in $b$. The restricted comeasuring bialgebra is
just $M_0=\C[t]/t^2=1$ and $\Delta t=t\tens t$, $\eps(t)=1$, and is
a Hopf algebra with $St=t$.

For $x^3=1$, the comeasuring bialgebra $M$ is generated by
$1,b,t,s$ modulo the relations
\[ b\alpha+t\gamma+s\beta=1,\quad t\alpha+s\gamma+b\beta=0,
\quad s\alpha+b\gamma+t\beta=0\]
where \[ \alpha\equiv b^2+ts+st,\quad
\beta\equiv s^2+bt+tb,\quad \gamma\equiv t^2+sb+bs.\]
The coalgebra is
\[ \Delta b=1\tens b+b\tens t+\alpha\tens s,\quad \Delta t=t\tens t
+\beta\tens s,\quad
\Delta s=s\tens t+\gamma\tens s,\quad \eps(b)=\eps(s)=0,\quad
\eps(t)=1,\]
i.e., the matrix of generators has the form
\[ \pmatrix{1&b&\alpha\cr 0&t&\beta\cr 0&s&\gamma}.\]
The extended $M_1$ is similar with weaker relations. Finally, the
restricted $M_0$ is generated by $1,t,s$ with the relations and
coproduct
\[ts+st=0,\quad t^3+s^3=1,\quad t^2s=ts^2=0,\quad
\Delta t=t\tens t+ s^2\tens s,\quad \Delta s=s\tens t+t^2\tens s\]
and $\eps(t)=1$, $\eps(s)=0$, i.e. the matrix of generators has the
form $\pmatrix{t&s^2\cr s &t^2}$. The second relation is equivalent
to $(t+s)^3=1$ given the others.

Similarly, one may compute the measuring bialgebra for
$k_\lambda=k[\lambda]/m(\lambda)=0$ for general fields and general
polynomials $m$. When $m$ is monic and irreducible, $k_\lambda$ is
a field extension of $k$ and $M(k_\lambda)$ should be viewed as the
`quantum Galois group' for the field extension. Also, we are not
limited to commutative algebras and fields. The comeasuring
bialgebras for the complex numbers and the quaternions, as algebras
of dimension 2, 4 respectively over $\R$ are computed in
\cite{AlbMa:qua} as preludes to the octonion case. In the
quaternion case one has a noncommutative and noncocommutative
bialgebra projecting onto the coordinate ring of $SO_3$. There are
of course plenty of other interesting examples according to ones
favourite algebra.

\subsection{Comeasurings of finite sets}

For completeness, we conclude our collection of general classes of
examples with the case $A=\C(\Sigma)$ where $\Sigma$ is a finite
set. In this case we compute only $M_1(\C(\Sigma))$ and
$M(\C(\Sigma))$ since there is no particularly natural splitting of
the identity without more structure. We take basis
$\{e_i=\delta_i,\ i\in
\Sigma\}$, the delta-functions at element in the set. The structure
constants are
\eqn{conset}{c_{ij}{}^k=\delta_{i,j}\delta^k_j}
and hence we find from Section~2 that $M_1(\C(\Sigma))$ is
generated by 1 and $\tau^i{}_j$ (say) with relations
\eqn{relset}{ \tau^k{}_i \tau^k_j=\delta_{i,j}\tau^k{}_i}
(no sum) and the matrix coalgebra structure. This means for each
row $k$ the matrix of generators forms an orthogonal family of
projectors, i.e. a copy of $\C(\Sigma)$, while there are no
relations between the different rows. We note that such algebras
have recently been considered in \cite{Wan:sym} in an interesting
$C^*$-algebra setting, again as some kind of universal automorphism
objects for finite sets. In our case we obtain them as an
elementary example of the general (but algebraic) construction in
Proposition~2.2 dual to \cite{Swe:hop}. We also have a unit
$1=\sum_i
\delta_i$ and the corresponding $M(\C(\Sigma))$ can be computed as
follows. We take a different basis $e_0=1$ and $\{e_i|\ i\in
\Sigma-*\}$ where $*$ is a basepoint (a fixed element of $\Sigma$).
Then $M(\C(\Sigma))$ is generated by $t^i{}_j,b_i$ with $i,j\in
\Sigma-*$ and relations
\eqn{reltbset}{ t^k{}_i t^k{}_j+b_it^k{}_j+t^k{}_ib_j
=\delta_{ij}t^k{}_j,\quad b_ib_j=\delta_{ij}b_j}
and the coalgebra in Proposition~2.3. We see that
$M(\C(\Sigma))\supset\C(\Sigma)$ embedded as $\{1,b_i\}$. Note that
if we want to describe $M(\C(\Sigma))$ in our original
basepoint-free delta-function basis, it consists of the quotient of
(\ref{relset}) by the relations
\eqn{quotset}{ \sum_j \tau^i{}_j=1} where $i,j\in\Sigma$. This is
equivalent to (\ref{reltbset}) via
\eqn{equivset}{\tau^i{}_j=t^i{}_j+b_j,\quad \tau^0{}_0=1-\sum_i b_i,
\quad \tau^0{}_i=b_i,\quad \tau^i{}_0=1-\sum_j t^i{}_j-\sum_j b_j}
for $i,j\in\Sigma-*$, which is the transformation induced by the
change between the two bases.

We note that the algebras $A=\C[x]/x^N=1$ are isomorphic by Fourier
transform to $\C(\Z_N)$ so their comeasuring bialgebras are
isomorphic to those for a finite set with $N$ elements. Explicitly,
\eqn{setZN}{ \tau^i{}_j=N^{-1}\sum_{m,n=0}^{N-1}
e^{{2\pi\imath\over N}(mi-nj)}t^m{}_n}
is the matrix of projectors (\ref{relset}) generating
$M_1(\C(\Sigma))$ in terms of the matrix generators of $M_1$ for
$\C[x]/x^N=1$ in the preceding subsection. For example, for $N=2$
the matrix of projectors is
\eqn{M1set2}{ \tau={1\over 2}\pmatrix{(b+t)^2+(b+t)&(b+t)^2-(b+t)\cr
(b-t)^2+(b-t)&(b-t)^2-(b-t)}.} That the elements of each row are
orthogonal projectors is equivalent to the relations $(b\pm
t)^3=(b\pm t)$. For the quotient $M(\C(\Sigma))$ we obtain
\eqn{Mset2}{ \tau=\pmatrix{g_+&1-g_+\cr 1-g_-&g_-},\quad g_\pm
={1\over 2}(1\pm (b\pm t))}
i.e. the algebra generated by two projectors $g_\pm$ and the
coalgebra
\eqn{Mset2delta}{ \Delta g_\pm=g_\pm\tens g_\pm+(1-g_\pm)
\tens (1-g_\mp),\quad \eps(g_\pm)=1}
as isomorphic to $M$ for $\C[x]/x^2=1$ in the preceding subsection.
This is the basepoint-free description via (\ref{relset}) and
(\ref{quotset}). On the other hand, if we compute $M(\C(\Sigma))$
from Proposition~2.3 according to (\ref{reltbset}) (which is in a
different basis) we have an equivalent matrix of generators
\[ \pmatrix{1&{}_\Sigma b\cr 0&{}_\Sigma t},
\quad {}_\Sigma b=1-g_+={1\over 2}(1-(b+t)),\quad {}_\Sigma t
=g_--1+g_+=t\]
obeying
\[ {}_\Sigma b^2={}_\Sigma b,\quad  ({}_\Sigma t+ {}_\Sigma b)^2
={}_\Sigma t+ {}_\Sigma b,
\quad \Delta {}_\Sigma b=1\tens {}_\Sigma b+{}_\Sigma b
\tens {}_\Sigma t,
\quad \Delta {}_\Sigma t={}_\Sigma t\tens {}_\Sigma t\] and
$\eps({}_\Sigma t)=1$,
$\eps({}_\Sigma b)=0$.

\subsection{Elements of quantum geometry}

In the remaining two subsections we consider briefly some steps
towards a `differential geometry' based on these `quantum
diffeomorphisms'. We essentially use the quantum group approach to
noncommutative geometry in \cite{BrzMa:gau} and especially the
recent paper \cite{Ma:rie}. We recall that in the classical
situation, if $\Sigma$ is a manifold  one has, roughly speaking, an
identification $\Diff(\Sigma)/\Diff_*(\Sigma)\isom
\Sigma$ providing a principal bundle over $\Sigma$ with structure
group $\Diff_*(\Sigma)$. The canonical projection sends a
diffeomorphism $\sigma$ to $\sigma(*)\in \Sigma$. Moreover, this
principal bundle has, again formally, a canonical form $\theta$
making this a frame resolution of $\Sigma$ in the language of
\cite{Ma:rie}. The usual affine frame bundle and linear frame
bundle are subbundles. Here $\theta$ is the projection to
$T_*\Sigma=\R^n$ of the Maurer-Cartan form on $\Diff(\Sigma)$,
where the Lie algebra of vector fields on $\Sigma$ modulo those
that vanish at $*$ is identified with the value at $*$. Using the
canonical form $\theta$, we have, essentially a correspondence
between gauge-fields on this principal bundle and covariant
derivatives on the cotangent bundle of $\Sigma$. In this way,
differential geometry on the manifold may be developed strictly as
$\Diff_*(\Sigma)$ gauge theory. Of course, because these are large
infinite dimensional groups, one needs topological considerations
to make these ideas fully precise in the manifold setting. One the
other hand, in our algebraic setting we can try to keep everything
algebraic, e.g if the role of our manifold is played by a
finite-dimensional (possibly noncommutative) algebra. We have seen
above that  for polynomials, $M(\C[x])$ plays the role of
`diffeomorphisms' and $M_0(\C[x])$ plays the role of
diffeomorphisms fixing the origin. Motivated by these ideas we make
the parallel choice for our finite dimensional algebras.

To recall the set-up, note  that arrows are reversed in our
co-ordinate ring formulation. Thus, a homogeneous quantum principal
bundle\cite{BrzMa:gau} arises from a Hopf algebra map $\pi: M\to
M_0$ between Hopf algebras $M,M_0$ (say) such that the induced
coaction $\Delta_L=(\pi\tens\id)\Delta:M\to M_0\tens M$ has fixed
point subalgebra $A=M^{M_0}\equiv\{h\in M|\ \Delta_L h=1\tens h\}$
(it plays the role of the coordinates of the base manifold), and
obeying a certain nondegeneracy condition (the `Hopf-Galois'
condition) that the map
\[ \chi:M\tens_A M\to M_0\tens M,\quad h\tens g\to \Delta_L(h)g\]
is invertible. There is a canonical form $\theta:V\to \Omega^1M$
given by\cite{Ma:rie}
\[ V=\ker\eps|_A,\quad \theta=(\id\tens S)\Delta \]
in terms of the Hopf algebra structure of $M$. Here
$\Omega^1M\subset M\tens M$ is the kernel of the product map (the
universal differential calculus). A connection on the quantum
bundle is a map $\omega:M_0\to\Omega^1M$ which is left
$\Ad$-covariant and obeys $\chi\circ\omega=\id\tens 1$ on
$\ker\eps\subset M_0$ and $\omega(1)=0$ as in \cite{BrzMa:gau}, and
in the presence of $\theta$ it defines a covariant derivative
$\nabla:\Omega^1A\to\Omega^1A\tens_A\Omega^1A\isom \Omega^2A\subset
A\tens A\tens A$ by
\[ \nabla\extd a=\extd a\tens 1 -a\o\tens Sa\t\omega(a\th)a\fo
-a\o\tens Sa\t\tens a\th+1\tens 1\tens a\]
where $\extd:A\to\Omega^1A$ given by $\extd a=a\tens 1-1\tens a$ is
the exterior derivative and $\Delta a=a\o\tens a\t$ etc., is a
notation for the coproduct, see \cite{Ma:rie}.

As a first step towards applying this formalism, we consider the
case where $\Sigma$ is a 2-point set. We take $M(\C(\Sigma))$ in
its basepoint form (\ref{reltbset}), namely generated by projectors
$p\equiv{}_\Sigma b+ {}_\Sigma t$ and $q\equiv {}_\Sigma b$ with no
further relations and the coalgebra
 \[ \Delta p=p\tens p+(1-p)\tens q,\quad\Delta
q=(1-q)\tens q+q\tens p,\quad \eps(p)=1,\quad
\eps(q)=1.\]
Similarly,  $M_0(\C(\Sigma))=\C[\bar p]/\bar p^2=\bar p$ with
$\Delta\bar p=\bar p\tens
\bar p$ and $\eps \bar p=1$. The projection is $\pi(p)=\bar p$ and
$\pi(q)=0$. The induced coaction is
\[ \Delta_Lq=1\tens q,\quad \Delta_Lp=\bar p\tens p
+(1-\bar p)\tens q.\]
Since $M$ is spanned by $1$ and alternating words $h=pqpq\cdots$ or
$h=qpqp\cdots$, one can show that
\[ \Delta_L1=1\tens 1,\quad \Delta_L h=\bar p\tens h
+(1-\bar p)\tens q\]
for either kind of non-trivial $h$. Hence the fixed point
subalgebra is the subalgebra spanned by $1,q$, i.e.
\[ M(\C(\Sigma))^{M_0(\C(\Sigma))}=\C(\Sigma).\]

One can also reach a similar conclusion
\[M(\C[x])^{M_0(\C[x])}=\C[x]\]
where $M(\C[x])=\C<t_i|\ i\in\Z_+>$ and $M_0(\C[x])=\C<\bar
t_i|i\in
\N>$ as explained in Section~3.2. The projection is $\pi(t_i)=\bar
t_i$ for $i>0$ and $\pi(t_0)=0$. The induced coaction is therefore
\[ \Delta_L t_0=1\tens t_0,\quad \Delta_L t_i
=\sum_{j=1}^{j=i}\sum_{n_1+\cdots +n_j=i}\bar t_{n_1}\cdots
\bar t_{n_j}\tens t_j\]
where all indices shown as in $\N$. Note that although the
coproduct of $M(\C[x])$ involves formal powerseries, its product
structure and the induced coaction are algebraic. Since $M(\C[x])$
and $M_0(\C[x])$ are free algebras, the fixed point subalgebra is
$\C[t_0]$.

To go further, one needs to have versions of $M,M_0$ which are
actually Hopf algebras and not merely bialgebras. This can
typically be done by localizing, i.e. by adjoining suitable
inverses or powerseries. For example, we may take the same $M$ as
for the two-point set but in the form of the comeasuring bialgebra
for $\C[x]/x^2=1$ as in Section~3.4 (generated by $1,b,t$). For the
splitting associated with this basis we have $M_0=\C[\bar t]/\bar
t^2=1$ which is actually a Hopf algebra. The map $\pi$ is
$\pi(b)=0$ and $\pi(t)=\bar t$. The induced coaction is
\[ \Delta_L b=1\tens b,\quad\Delta_L t=\bar t\tens t\]
Since $t^2=1-b^2$ in $M$, every element of that can be written in
the form $f(b)+tg(b)$, of which the fixed elements are those with
$g=0$. Therefore the fixed subalgebra $A$, the base of the bundle,
is
\[  M^{M_0}=\C[b]\]
Note that the base is no longer the original algebra of which we
took the `diffeomorphisms'. (This is attributable to the nontrivial
$c_{ij}{}^0$ in Proposition~2.3 for this basis). On the other hand,
if we allow formal powerseries in $b^2$ then $M$ has an antipode
and, moreover, the nondegeneracy (Hopf-Galois) condition holds so
that we have, at least formally, a homogeneous quantum principal
bundle. The canonical form $\theta$ is also formal, for the same
reason. Here $\ker\eps=\{f(b)+tg(b)|\ f(0)+g(0)=0\}$ and hence $V$
consists of polynomials in $b$ vanishing at 0. Then
\[ \theta(b)=b\tens {t\over 1-b^2}-1\tens {bt\over 1-b^2}.\]
On the other hand, the map $\Phi(\bar t)=t$ is a coalgebra map
splitting of $\pi$ and makes this localized bundle trivial. It
should therefore be viewed only as a local `patch' of a more non
trivial bundle without such localization. (This requires, however,
a suitable extension of the theory in
\cite{BrzMa:gau}\cite{Ma:rie}.) At least in this patch, a
connection $\omega$ is induced as in \cite{BrzMa:gau} by a `gauge
field' $\alpha:M_0\to \Omega^1\C[b]$ such that $\alpha(1)=0$, which
in our case means a single element $\alpha=\alpha(\bar t)\in
\Omega^1\C[b]$. This induces the
covariant derivative  \[
\nabla:\Omega^1\C[b]\to\Omega^1\C[b]\tens_{\C[b]}\Omega^1\C[b],\quad
 \nabla\extd b=(b\tens 1-1\tens b)\alpha.\]

\subsection{Preserving nonuniversal differential calculi}

Here we discuss briefly the further restrictions imposed by
non-universal differential calculi. We recall that for any unital
algebra $A$ the universal calculus is $\Omega^1A$ is the kernel of
the product map $A\tens A\to A$ as an $A$-bimodule, while a general
differential calculus is a quotient $\Omega^1(A)=\Omega^1A/N$ for
some subbimodule $N$. The differential is the universal one $\extd
a=a\tens 1-1\tens a$ (as above) projected down to $\Omega^1(A)$ in
the nonuniversal case. It is natural to restrict the notion of
unital comeasurings to ones that preserve $N$ in the sense
\eqn{difcomeas}{\beta(N)\subseteq N\tens B,}
where $\beta$ is extended to the tensor square $A\tens A$ in the
usual way,  restricted to $\Omega^1A$. The universal object in this
case is a quotient of $M(A)$. It is easy to see that it remains a
bialgebra, which we denote $M(A,\Omega^1(A))$, the {\em comeasuring
bialgebra with nonuniversal calculus}.

Thus, when $A=\C(\Sigma)$, the subbimodules $N$ are classified are
classified by subsets of $\Sigma\times\Sigma-{\rm diag}$. In
\cite{BrzMa:dif} the complement of this subset is denoted $E$ and
we write $i-j$ when $(i,j)\in E$ and $i\# j$ when $(i,j)$ is in the
complement in $\Sigma\times\Sigma-{\rm diag}$. Thus,
$N=\{\delta_i\tens\delta_j|\ i\# j\}$. To describe the comeasuring
bialgebra $M(\C(\Sigma),\Omega^1(\C(\Sigma)))$ in this case we use
the base-point free version of $M(\C(\Sigma))$ defined by
(\ref{relset}) and (\ref{quotset}), and quotient further by the
relations
\eqn{difset}{ \tau^i{}_j \tau^k{}_l=0\quad \forall i-k,\quad j\#l}
as the requirement that $N$ is preserved by $\delta_j\mapsto
\delta_i\tens \tau^i{}_j$. Note that in $M(\C(\Sigma))$ we have,
for $i-k$ and $j\#l$,
\[ \Delta (\tau^i{}_j \tau^k{}_l)=\sum_{a-b} \tau^i{}_a \tau^k{}_b\tens
\tau^a{}_j \tau^b{}_l
+\sum_{a\# b} \tau^i{}_a \tau^k{}_b\tens \tau^a{}_j \tau^b{}_l.\]
Hence it is clear that the quotient by the relations (\ref{difset})
remains a bialgebra. The case $a=b$ does not contribute here in
view of $j\ne l$ and (\ref{relset}).

When $\Sigma$ is the two-point set, there  is (up to equivalence)
only one nonuniversal nontrivial differential calculus, namely
$E=\{(1,2)\}$. Then $M(\C(\Sigma),\Omega^1(\C(\Sigma)))$ is the
quotient of the bialgebra (\ref{Mset2}), (\ref{Mset2delta}) by the
additional relation
\[ (1-g_+)(1-g_-)=0.\]
Equivalently, we quotient the comeasuring bialgebra $M$ for
$\C[x]/x^2=1$ in Section~3.4 by
\[ t^2=t(1+b).\]

When $A=\C[x]$ the natural differential calculi are those
bicovariant under the coaddition structure, and are labeled by a
single parameter $\lambda\in\C$ (and over a general field, the
coirreducible calculi correspond to field extensions). The standard
commutative differential calculus is the one where $N$ is the
subbimodule generated by $x\extd x-(\extd x)x$. This contains in
particular all elements of the form $x^i\extd x^j-(\extd x^j)x^i$
for $i,j\in
\Z_+$. Under the coaction (\ref{coactpoly}) we have
\[ x\extd x-(\extd x)x\mapsto \sum_{i,j}(x^i\extd x^j
-(\extd x^j)x^i)\tens t_it_j+
\sum (\extd x^j)x^i\tens (t_it_j-t_jt_i).\]
Since the $(\extd x^j)x^i$ are linearly independent elements of
$\Omega^1\C[x]$ for $j>0$, the differentiable comeasurings
$M(\C[x],\Omega^1(\C[x]))$ for the standard commutative
differential calculi are precisely the quotient of
$M(\C[x])=\C<t_i|\ i\in\Z_+>$ by the relations of commutativity of
the generators i.e. precisely the usual $\Diff(\C[x])$.
Differential calculi in between the universal one and the
commutative one therefore lead to quotients in between the free
comeasuring bialgebra $M(\C[x])$ and the classical commutative one.

\subsection{Diffeomorphisms of the quantum-braided plane}

Until now, we have only considered diffeomorphism of the line and
its quotients. We now move on to the plane, which is now nontrivial
enough to have a noncommutative $q$-deformed version, the
quantum-braided plane $\C_q^2$. This has generators $1,x,y$ with
relations $yx=qxy$, where $q\ne0$ is a parameter. It is
infinite-dimensional but, as in Section~3.1, we find that
$M_0(\C_q^2)$ is algebraic without the need for formal powerseries.

The computation of the comeasuring bialgebra $M(\C[x,y])$ for the
classical plane follows just the same steps as in Section~3.1. With
the basis $\{e_{m,n}=x^my^n|\ m,n\in\Z_+\}$ one has the structure
constants
\[ c_{(i,j)(k,l)}{}^{(m,n)}=\delta^m_{i+k}\delta^n_{j+l}\]
and $M_1(\C[x,y])$ is generated by $1,t^{(i,j)}{}_{(k,l)}$ with
relations
\eqn{xyrel}{ \sum_{(a,b)+(c,d)=(m,n)} t^{(a,b)}{}_{(i,j)}
t^{(c,d)}{}_{(k,l)}
=t^{(m,n)}{}_{(i,j)+(k,l)}}
and the matrix coalgebra. Similarly to Section~3.1, we find
$M(\C[x,y])$ is generated by $1$, $s_{(i,j)},t_{(i,j)}$ where
$i,j\in\Z_+$ modulo the relations
\eqn{strel}{ \sum_{(a,b)+(c,d)=(i,j)} s_{(a,b)}t_{(c,d)}
= \sum_{(a,b)+(c,d)=(i,j)}  t_{(a,b)}s_{(c,d)},\quad
\forall i,j\in\Z_+.}
Here
\eqn{stdef}{ s_{(i,j)}\equiv t^{(i,j)}{}_{(1,0)},\quad
t_{(i,j)}\equiv t^{(i,j)}{}_{(0,1)}}
generate the others by repeated `convolution'. Thus,
\[ t^{(i,j)}{}_{(k,0)}=\sum_{(a_1,b_1)+\cdots(a_k,b_k)
=(i,j)} s_{(a_1,b_1)}
\cdots s_{(a_k,b_k)},\quad
t^{(i,j)}{}_{(0,k)}=\sum_{(a_1,b_1)+\cdots(a_k,b_k)=(i,j)}
t_{(a_1,b_1)}\cdots t_{(a_k,b_k)}\]
\[ t^{(i,j)}{}_{(k,l)}=
\sum_{(a,b)+(c,d)=(i,j)} t^{(a,b)}{}_{(k,0)}t^{(c,d)}{}_{(0,l)}\]
for $k,l\ge 1$. These follow from (\ref{xyrel}), while
(\ref{strel}) is the residual content of (\ref{xyrel}) after making
these substitutions. If we consider the $s_{(i,j)}$, $t_{(i,j)}$ as
sequences on $\Z_+\times\Z_+$ with convolution product
$s*t_{(i,j)}=\sum_{(a,b)+(c,d)=(i,j)}s_{(a,b)} t_{(c,d)}$, we can
write the residual relations of $M(\C[x,y])$ compactly as
$s*t=t*s$. The coalgebra is
\eqn{stdelta}{\Delta s_{(i,j)}=\sum_{(a,b)}
t^{(i,j)}{}_{(a,b)}\tens s_{(a,b)},\quad
\Delta t_{(i,j)}=\sum_{(a,b)} t^{(i,j)}{}_{(a,b)}\tens t_{(a,b)},
\quad \eps (s_{(i,j)})=\delta^i_1\delta^j_0,\quad
 \eps (t_{(i,j)})=\delta^i_0\delta^j_1.}
The coaction on $\C[x,y]$ is
\eqn{stcoact}{ x\mapsto \sum_{(i,j)} x^iy^j\tens s_{(i,j)},
\quad y\mapsto
\sum_{(i,j)} x^iy^j\tens t_{(i,j)}.}
The quotient $M_0(\C[x,y])$ has the same form but with the indices
in (\ref{stdef}), (\ref{stdelta}) and (\ref{stcoact}) excluding
$(0,0)$. In this case $t^{(i,j)}{}_{(k,l)}=0$ unless $i+j\ge k+l$,
so that the coproduct in this case is a finite sum. The geometric
meaning of $M_0(\C[x,y])$ is the diffeomorphisms that fix the
origin.

For the quantum-braided plane, we again have a basis
$\{e_{m,n}=x^my^n|\ m,n\in\Z_+\}$, but now with the q-deformed
structure constants and consequent relations
\eqn{qxyrel}{c_{(i,j)(k,l)}{}^{(m,n)}=\delta^m_{i+k}
\delta^n_{j+l}q^{jk},\quad
\sum_{(a,b)+(c,d)=(m,n)}q^{bc} t^{(a,b)}{}_{(i,j)}t^{(c,d)}{}_{(k,l)}
=q^{jk}t^{(m,n)}{}_{(i,j)+(k,l)}.}

\begin{propos} The comeasuring bialgebra $M(\C_q^2)$ of the
quantum-braided plane
is generated by $1$, $s_{(i,j)},t_{(i,j)}$ for $i,j\in\Z_+$, modulo
the relations
\[
q\sum_{(a,b)+(c,d)=(i,j)} q^{bc}s_{(a,b)}t_{(c,d)}
= \sum_{(a,b)+(c,d)=(i,j)} q^{bc} t_{(a,b)}s_{(c,d)},\quad
\forall i,j\in\Z_+.\]
The coaction is as in (\ref{stcoact}). The higher generators and
hence the coproduct (\ref{stdelta}) are given by
\[ t^{(i,j)}{}_{(k,0)}=\sum_{(a_1,b_1)+\cdots(a_k,b_k)=(i,j)}
s_{(a_1,b_1)}\cdots s_{(a_k,b_k)} q^{\sum_{s=2}^k(b_1+\cdots
+b_{s-1})a_s}\]\[
t^{(i,j)}{}_{(0,k)}=\sum_{(a_1,b_1)+\cdots(a_k,b_k)=(i,j)}
t_{(a_1,b_1)}\cdots t_{(a_k,b_k)} q^{\sum_{s=2}^k(b_1+\cdots
+b_{s-1})a_s}\]
\[ t^{(i,j)}{}_{(k,l)}=\sum_{(a,b)+(c,d)=(i,j)} q^{bc}
t^{(a,b)}{}_{(k,0)}t^{(c,d)}{}_{(0,l)}
\] for $k\ge 1$.
\end{propos}
\proof The relations (\ref{qxyrel}) allow one to define the
general $t^{(i,j)}{}_{(k,l)}$ as stated in terms
of the $s_{(i,j)}$, $t_{(i,j)}$. Putting these back into
(\ref{qxyrel}) leaves the residual relations between $s,t$ as
shown. These are obtained from
$t^{(i,j)}_{(1,0)+(0,1)}=t^{(i,j)}_{(1,1)}=t^{(i,j)}_{(0,1)+(1,0)}$
computed from (\ref{qxyrel}), i.e. they reflect commutativity of
the addition law on $\Z_+\times\Z_+$. To see that these are all the
relations, it is useful to define the $q$-deformed convolution
product $s*_qt_{(i,j)}=
\sum_{(a,b)+(c,d)=(i,j)} q^{bc}s_{(a,b)}t_{(c,d)}$ for sequences
on $\Z_+\times\Z_+$.
One may check that $*_q$ is associative. Then
$t^{(i,j)}{}_{(k,0)}=(s*_q\cdots *_qs)_{(i,j)}$ and
$t^{(i,j)}{}_{(0,k)}=(t*_q\cdots *_qt)_{(i,j)}$ ($k$-fold products)
and $t^{(i,j)}{}_{(k,l)}=(s*_q\cdots *_q s*_q t*q\cdots *_q
t)_{(i,j)}$ ($k$-fold and $l$-fold). The relations (\ref{qxyrel})
then hold if $t*_q s=qs*_qt$, which are the stated relations of
$M(\C_q^2)$. \endproof

The quotient $M_0(\C_q^2)$ has the same form with the index $(0,0)$
excluded  from the expressions in the proposition. In this case
$t^{(i,j)}{}_{(k,l)}=0$ unless $i+j\ge k+l$, so that the coproduct
is a finite sum. The lowest level generators of $M_0(\C_q^2)$ are
\eqn{qmatst}{ \pmatrix{a&b\cr c&d}\equiv\pmatrix{s_{(1,0)} &
t_{(1,0)}\cr s_{(0,1)}& t_{(0,1)}}
=\pmatrix{t^{(1,0)}{}_{(1,0)}& t^{(1,0)}{}_{(0,1)}\cr
t^{(0,1)}{}_{(1,0)}&t^{(0,1)}{}_{(0,1)}}} and the relations among
these are
\[ dc=qcd,\quad ba=qab,\quad ad-da=q^{-1}bc-qcb,\]
which are just half of the relations of the quantum matrices.
Moreover, this is just the lowest level content of the relations
$t*_qs=qs*_qt$.

\subsection{Preserving a coaddition}

We are now ready to consider the general theory of restricting
diffeomorphisms to those preserving a coalgebra structure. In our
`coordinate ring' setting it means a quotient of $M_1(A)$, which we
denote $M_1(A,\Delta)$. Thus, a coalgebra structure on $A$
corresponds in a basis $\{ e_i\}$ to structure constants defined by
$\Delta e_i=d_i{}^{jk}e_j\tens e_k$. Along the same lines as in
Section~2, the comeasurings that preserve this clearly means the
quotient of $M_1(A)$ by the additional relation
\eqn{tdrel}{ d_a{}^{jk} t^a{}_i=d_i{}^{ab} t^j{}_a t^k{}_b.}
As explained below Proposition~2.2,  one has a bialgebra for any
$d_i{}^{jk}$, i.e. one does not need $\Delta$ to be coassociative
or to make $A$ into a bialgebra. The coassociativity is, however,
natural to assume.

We also have quotients $M(A,\Delta),M_0(A,\Delta)$ etc., as before.
And in principle, we have still further quotients where a counit
$\eps$ is respected as well. However, the counit in the case of a
bialgebra $A$ defines a natural splitting $A=1\oplus A'$, where
$A'=\ker\eps$, i.e. it is natural to choose the basis so that
$\eps(1)=1$ and $\eps(e_i)=0$ for $i>0$. At least in this case, the
bialgebra $M_0(A,\Delta)$ automatically preserves the counit
without further quotients arising from that.

The quantum-braided plane $\C_q^2$ has such coalgebra structure,
expressing the braided addition law. This is the braided
`coaddition' introduced in \cite{Ma:poi}. Explicitly, it is given
by
\[ \Delta_+ (x^m y^n)=\sum_{r=0}^m\sum_{s=0}^n\left[{m\atop r}
\right]_{q^2}x^ry^s\tens
x^{m-r}y^{n-s} q^{(m-r)s},\] where we use the $q^2$-binomial
coefficients defined as usual, but in terms of $q^2$-factorials
\[ [n]_{q^2}!=[n]_{q^2}\cdots [1]_{q^2},\quad [m]_{q^2}={1-q^{2m}
\over 1-q^2}.\]
Although this coproduct forms a coalgebra (it is coassociative), it
does not make $\C_q^2$ into a usual bialgebra or Hopf algebra, but
rather into a braided group\cite{Ma:bra}. We will say more about
this in the next section; for the present purposes we need only to
know its explicit form as stated here. For our above basis, we have
\[ d_{(m,n)}{}^{(i,j)(k,l)}=\delta^{i+k}_m\delta^{j+l}_n
\left[{m\atop i}\right]_{q^2}
\left[{n\atop j}\right]_{q^2}q^{jk}.\]

\begin{propos} The quotient of the restricted comeasuring bialgebra
$M_0(\C_q^2)$
respecting the braided coaddition on $\C_q^2$, for generic $q$, can
be identified with the standard $2\times 2$ quantum matrices
$M_q(2)$.
\end{propos}
\proof  The additional quotient of $M_0(\C_q^2)$ is by the relation
\[\sum_{(a,b)+(c,d)=(m,n)} \left[{m\atop a}\right]_{q^2}
\left[{n\atop b}\right]_{q^2}q^{bc}t^{(i,j)}{}_{(a,b)}
t^{(k,l)}{}_{(c,d)}=
\left[{i+k\atop i}\right]_{q^2}
\left[{j+l\atop j}\right]_{q^2}q^{jk} t^{(i,j)+(k,l)}{}_{(m,n)}.\]
For generic $q$ we can define new generators
\[\tau^{(i,j)}{}_{(k,l)}=t^{(i,j)}{}_{(k,l)}{[i]_{q^2}!
[j]_{q^2}!\over[k]_{q^2}![l]_{q^2}!}\]
and then the additional relations become
\[\sum_{(a,b)+(c,d)=(m,n)}q^{bc} \tau^{(i,j)}{}_{(a,b)}
\tau^{(k,l)}{}_{(c,d)}
=q^{jk}\tau^{(i,j)+(k,l)}{}_{(m,n)}\]
which, by similar reasoning as for the comeasuring bialgebra,
implies that the generators can be obtained by convolution from the
generators $\sigma_{(i,j)}=\tau^{(1,0)}{}_{(i,j)}$ and
$\tau_{(i,j)}=\tau^{(0,1)}{}_{(i,j)}$. This time
\[ t^{(k,0)}{}_{(i,j)}=\sum_{(a_1,b_1)+\cdots(a_k,b_k)=(i,j)}
\sigma_{(a_1,b_1)}\cdots \sigma_{(a_k,b_k)} q^{\sum_{s=2}^k(b_1+\cdots
+b_{s-1})a_s}\]\[
t^{(0,k)}{}_{(i,j)}=\sum_{(a_1,b_1)+\cdots(a_k,b_k)=(i,j)}
\tau_{(a_1,b_1)}\cdots \tau_{(a_k,b_k)} q^{\sum_{s=2}^k(b_1+\cdots
+b_{s-1})a_s} \]
\[ \tau^{(k,l)}{}_{(i,j)}=\sum_{(a,b)+(c,d)=(i,j)} q^{bc}
\tau^{(k,0)}{}_{(a,b)}t^{(0,l)}{}_{(c,d)}
\] for $k,l\ge 1$. The residual relations are
\[
q\sum_{(a,b)+(c,d)=(i,j)} q^{bc}\sigma_{(a,b)}\tau_{(c,d)}
= \sum_{(a,b)+(c,d)=(m,n)} q^{bc} \tau_{(a,b)}\sigma_{(c,d)},
\quad \forall i,j\in\Z_+\]
of in the convolution notation, $\tau*_q\sigma=q\sigma*_q\tau$.
This is the bialgebra respecting only the coalgebra $\Delta_+$
(indeed, the quantum-braided plane is self-dual as a braided group
and this is why this algebra has the same form as the comeasuring
bialgebra).

These convolution formulae imply that $t^{(i,j)}{}_{(k,l)}=0$
unless $i+j\le k+l$. Combined with the reverse inequality  for
$M_0(\C_q^2)$, we see that $t^{(i,j)}{}_{(k,l)}=0$ in
$M(\C_q^2,\Delta_+)$ unless $i+j=k+l$. It follows that
$M_0(\C_q^2,\Delta_+)$ is generated by $1$ and the lowest level
generators (\ref{qmatst}). Half their relations are given above,
inherited from $M_0(\C_q^2)$. The other half come from the
relations $\tau*_q\sigma=q\sigma*_q\tau$, which are computed
similarly as
\[ db=qbd,\quad ca=qac,\quad ad-da=q^{-1}cb-qbc.\]
Thus, $M_0(\C_q^2,\Delta_+)=M_q(2)$, the standard $2\times 2$
quantum matrices in the conventions of \cite{Ma:book}. The coaction
reduces to the standard coaction on $\C_q^2$. \endproof

\section{R-matrix constructions for comeasuring bialgebras}

In this section we consider some general constructions possible
when our algebra is braided, i.e. in the presence of a Yang-Baxter
operator. As a first application, we note that until now we have
studied the maximal comeasuring objects in the category of
bialgebras. Since they are maximal they tend to be free, with only
the minimal relations compatible with coacting on our algebra.
However, when the algebra is itself braided, we can look for
objects maximal in some braided-commutative sense. We then give
more functorial braided group versions of these constructions,
related by transmutation.

We recall that an algebra $A$ is braided if it comes with an
operator $\Psi:A\tens A\to A\tens A$ obeying the braid relations
and functorial with respect to the product of $A$ in the sense
\eqn{braalgpsi}{ \Psi(ab\tens c)=(\id\tens\cdot)(\Psi\tens\id)
(\id\tens\Psi)(a\tens b\tens c),
\quad \Psi(a\tens bc)=(\cdot\tens\id)(\id\tens\Psi)
(\Psi\tens\id)(a\tens b\tens c),}
for all $a,b,c\in A$. More generally, a braided algebra means $A$
an object in a braided category (such as that generated by a single
braiding operator) with the product a morphism. We recall that in a
braided category there is a braiding between any two objects
playing the role of `transposition'. Also, the (co)modules of any
(dual) quasitriangular bialgebra or Hopf algebra (i.e., of any
`strict' quantum group) form a braided category, so any algebra
covariant under a strict quantum group is a braided algebra. See
\cite{Ma:book} or papers such as \cite{Ma:bra}\cite{Ma:exa} where
braided algebras and groups have been introduced. In this setting
we introduce $M_1(R,A)$ as again a dual quasitriangular bialgebra
or `strict' quantum group.

Thus, when $A$ has a basis $\{e_i\}$, a braided algebra translates
into the existence of a matrix $R\in M_n\tens M_n$ (where $M_n$
denotes $n\times n$-matrices and $n=\dim(A)$) obeying the Quantum
Yang-Baxter Equations (QYBE)
$R_{12}R_{13}R_{23}=R_{23}R_{13}R_{12}$ (a so-called R-matrix),
such that
\eqn{braalgmat}{ c_{12}{}^3R_{14}R_{24}=R_{34}c_{12}{}^3,
\quad R_{12}c_{34}{}^2c_{34}{}^2R_{14}R_{13}}
We use the standard compact notation where $R_{12}=R\tens\id$ and
we suppose that $R$ is invertible. In explicit component terms, the
requirement is
\eqn{braalg}{ c_{ij}{}^a R^k{}_a{}^m{}_n=c_{ab}{}^k
R^a{}_i{}^m{}_c R^b{}_j{}^c{}_n,
\quad c_{jk}{}^a R^m{}_n{}^i{}_a=c_{ab}{}^i R^m{}_c{}^b{}_k
R^c{}_n{}^a{}_j,}
where
\[ \Psi(e_i\tens e_j)=e_b\tens e_a R^a{}_i{}^b{}_j.\]

\subsection{Dual quasitriangular comeasuring bialgebras}

We also recall that a bialgebra $M$ is dual quasitriangular or
$\CR$-commutative if there exists $\CR:M\tens M\to k$ such that
\[ \CR(ab,c)=\CR(a,c\o)\CR(b,c\t),\quad \CR(a,bc)=\CR(a\o,c)\CR(a\t,b),
\quad b\o a\o\CR(a\t,b\t)=\CR(a\o,b\o)a\t b\t\]
for all $a,b,c\in M$. The axioms are dual to the quasitriangular
structures introduced for the quantum groups $U_q(g)$ by
Drinfeld\cite{Dri:qua}. On the other hand, given any $R$-matrix
there is a dual quasitriangular bialgebra $A(R)$ of `quantum
matrices' with generators $R\vect_1\vect_2=\vect_2\vect_1R$ (in a
compact notation where $\vect_1=\vect\tens\id$ etc.). These are the
FRT relations\cite{FRT:lie} while the dual quasitriangularity for
general $R$ is due to the author\cite{Ma:qua} and takes the form
$\CR(t^i{}_j,t^k{}_l)=R^i{}_j{}^k{}_l$. See \cite{Ma:book}.

\begin{propos} If $A$ is a braided algebra with braiding defined by
an $R$-matrix, then $M_1(R,A)$ defined by the matrix coalgebra and
the relations
\[ R^i{}_a{}^j{}_b t^a{}_k t^b{}_l=t^j{}_b t^i{}_a R^a{}_k{}^b{}_l,
\quad c_{ij}{}^a t^k{}_a=c_{ab}{}^k t^a{}_i t^b{}_j\] is a dual
quasitriangular bialgebra.
\end{propos}
\proof  We check that this quotient of $A(R)$, in which the relations
of Proposition~2.2
are further imposed, inherits the linear functional $\CR$. If so
then it will remain a dual quasitriangular structure for the
quotient. Thus,
\cmath{\CR(t^i{}_j,c_{ab}{}^k t^a{}_m t^b{}_n)=R^i{}_c{}^b{}_n
R^c{}_j{}^a{}_m c_{ab}{}^k
=c_{mn}{}^aR^i{}_j{}^k{}_a=\CR(t^i{}_j,c_{mn}{}^at^k{}_a)\\
\CR(c_{ab}{}^kt^a{}_i t^b{}_j,t^m{}_n)=c_{ab}^kR^a{}_i{}^m{}_c
R^b{}_j{}^c{}_n=c_{ij}{}^a
R^k{}_a{}^m{}_n=\CR(c_{ij}{}^at^k{}_a,t^m{}_n)} using the
covariance conditions (\ref{braalg}). The proof for $t^i{}_j$ and
$t^m{}_n$ replaced by general strings of generators has just the
same form, with repeated use of the covariance conditions; the
general proof thereby proceeds by a straightforward induction.
\endproof

While this quotient will always be dual quasitriangular, it may
also be trivial, i.e. the ideal generated by both sets of relations
may be too large. Although not necessary, if $A$ is itself `braided
commutative' in some sense then one may expect that the above
construction is more natural. The appropriate form of commutativity
is, as for braided matrices and braided
planes\cite{Ma:exa}\cite{Ma:poi}, the existence of a matrix $R'\in
M_n\tens M_n$ obeying
\eqn{R'R}{ R'_{12}R_{13}R_{23}=R_{23}R_{13}R'_{12},
\quad R_{12}R_{13}R'_{23}=R'_{23}R_{13}R_{12}.}
More precisely, these matrix conditions follow from and are
essentially equivalent to the algebraic condition
$R'\vect_1\vect_2=\vect_2\vect_1R'$ in the quantum matrix algebra.
Typically, $R'$ is built from $R$ and the relations with $R'$ are
equivalent to the defining relations with $R$. With respect to
this, braided-commutativity is
\eqn{R'comm}{ c_{ij}{}^k=c_{ba}{}^k R'^a{}_i{}^b{}_j.}
Note that the more naive definition of braided commutativity would
be $\cdot\circ\Psi(a\tens b)=ab$ for all $a,b$, but this is only
natural when $\Psi^2=\id$. It corresponds to the choice $R'=R$,
which is too restrictive to apply in most examples of
q-deformation. Then
\cmath{ c_{ba}{}^k(t^b{}_jt^a{}_iR'^i{}_m{}^j{}_n)
=c_{ba}{}^kR'{}^a{}_i{}^b{}_jt^i{}_mt^j{}_n
=c_{ij}{}^kt^i{}_m t^j{}_n=c_{mn}{}^at^k{}_a\\
=c_{ji}{}^at^k{}_aR'{}^i{}_m{}^j{}_n
=(c_{ba}{}^k t^b{}_jt^a{}_i)R'^i{}_m{}^j{}_n}
holds automatically in $M_1(R,A)$.

As before, there are quotients $M(R,A)$ and $M_0(R,A)$ when $\Psi$
respects $e_0=1$. The natural condition is $\Psi(a\tens 1)=1\tens
a$ and $\Psi(1\tens a)=a\tens 1$ for all $a\in A$.

As a finite-dimensional example, the $q$-epsilon tensor and
$q$-metric associated to the $SO_q(3)$-covariant quantum plane may
be used to define an associative $q$-quaternion algebra $\H_q$. Its
automorphism bialgebra then recovers $M_0(R,\H_q)=SO_q(3)$ along
the same lines as for the $q=1$ case in \cite{AlbMa:qua}.

For the simplest `geometrical' example, namely $A=\C[x]$, we have
more than one way to consider it as a braided algebra. The simplest
is as the braided line\cite{Ma:poi}\cite{Ma:sol} where
\eqn{braline}{\Psi(x^i\tens x^j)=q^{ij}x^j\tens x^i,\quad
R^i{}_j{}^k{}_l
=\delta^i{}_j\delta^k{}_l q^{ik}.}
Then $M(R,\C[x])=\Diff(\C[x])$ is the classical commutative
diffeomorphism group. Here $q$ cancels from both sides and the
result is the same as for $q=1$, i.e. the imposition of
commutativity between all the generators.

A different braiding on $A=\C[x]$ is the one introduced in
\cite{Ma:sol},
\eqn{dbraline}{ \Psi(x^i\tens x^j)=\sum_{k=0}^i\left[{i\atop k}
\right]_qq^{j(i-k)}(1-q)^k{[j+k-1]_q!\over
[j-1]_q!} x^{j+k}\tens x^{i-k},\quad\forall i,j\in\N} and the
trivial transposition when $i$ or $j=0$. This is the canonical
`double' braiding associated to any braided group, in this case the
braided line. The algebra of the braided group is automatically a
braided algebra under its canonical braiding. From another point of
view, the double bosonisation of a braided group canonically acts
on the braided group. Here the double bosonisation of the braided
line is $U_q(su_2)$ and acts as a q-deformation of $so(2,1)$ by
`conformal transformation' on $\C[x]$\cite{Ma:conf}. The above
braiding is the one induced from the quasitriangular structure of
$U_q(su_2)$ by this action. We therefore call this the `conformally
braided' line.

\begin{propos} The restricted dual-quasitriangular comeasuring bialgebra
$M_0(R,\C[x])$ is generated by
$1,t_i$ for $i\in\N$ with the relations
\[ qt_jt_i=\sum_{k=0}^{j-1}{[i+k]_q![j-1]_q!\over [j-k-1]_q![i]_q!
[k]_q!}q^{i(j-k)}(1-q)^k t_{i+k}t_{j-k}\]
and the coalgebra (\ref{deltapoly}) as in Section~3.1. The dual
quasitriangular structure is $\CR(t_i,t_j)=q\delta^i_1\delta^j_1$.
\end{propos}
\proof The R-matrix corresponding to (\ref{dbraline}) is
\[
R^i{}_0{}^j{}_k=\delta^i{}_0\delta^j{}_k,\quad
R^i{}_j{}^k{}_0=\delta^i{}_j\delta^k{}_0,\quad
R^i{}_j{}^k{}_l=\delta^{i+k}_{j+l}\left[{j\atop
j-i}\right]_qq^{li}(1-q)^{j-i}{[k-1]_q!\over [l-1]_q!}\] for
$j,k\in \N$ and $i\le j$ in the third expression (which is zero
otherwise), which implies that
$R^a{}_1{}^b{}_1=q\delta^a_1\delta^b_1$ when $a,b\in\N$. Hence the
$R\vect_1\vect_2=\vect_2
\vect_1R$ relations on the generators $t_i\equiv t^i{}_1$ have the form
\[ \sum_{a,b}R^i{}_a{}^j{}_b t_a t_b=q t_j t_i,\]
which then computes as stated. The R-matrix also provides the
dual-quasitriangular structure on the generators.
\endproof

In particular, the relations of the lowest order generators with
the general generators are
\cmath{ t_1t_i=q^{i-1}t_i t_1,\quad t_2t_i=q^{2i-1}t_it_2
+(1-q^{i+1})q^{i-1}t_{i+1}t_1\\
t_3t_i=q^{3i-1}t_it_3+(1-q^{i+1})(1+q)q^{2i-1}t_{i+1}t_2+(1-q^{i+2})
(1-q^{i+1})q^{i-1}t_{i+2}t_1}
for all $i\in\N$, etc. From this, for generic $q$, one obtains
explicitly
\cmath{ t_1t_2=qt_2t_1,\quad t_1t_3=q^2t_3t_1=q(t_2)^2,
\quad t_2t_3=qt_3t_2=q^2t_4t_1,\\
\Delta t_1=t_1\tens t_1,\quad \Delta t_2=t_2\tens t_1+(t_1)^2\tens t_2,
\quad\Delta t_3=t_3\tens t_1+(1+q)t_2t_1\tens t_2+(t_1)^3\tens t_3}
and so on, to all orders. The comeasuring bialgebra $M(R,\C[x])$ is
similar, with the extra generator $t_0$. The coproduct in this case
involves infinite sums, so has to be treated formally.

For anyonic variables $x^3=0$, we again have a braided group (the
anyonic line), and hence a canonical `double' braiding. Here
$q^3=1$ and the resulting $9\times 9$ $R$-matrix in this case is
given in \cite[Ex.~4.8]{Ma:sol}. The restricted $M_0(R,\C[x]/x^3)$
is the quotient of $M_0(\C[x]/x^3)$ in Section~3.2 by the relation
$ts=qst$, i.e. we have a dual-quasitriangular bialgebra
\[ ts=qst,\quad \Delta t=t\tens t,\quad \Delta s=s\tens t+t^2\tens s,
\quad \CR(t^i,t^j)=q^{ij}\]
and $\CR$ zero on any expression containing $s$.

Finally, we let $A=\C_q^2$ be the quantum-braided plane $\C_q^2$
with generators $x,y$ and relations $yx=qxy$, as in Section~3.7.
This has a natural braiding $\Psi$ used to describe coaddition on
the quantum plane\cite{Ma:poi}, as studied in Section~3.8. The
quantum-braided plane is also braided commutative with respect to a
certain matrix $R'$. Explicitly,
\eqn{xybra}{ \Psi(x\tens x)=q^2x\tens x,\quad \Psi(y\tens y)
=q^2y\tens y,\quad\Psi(x\tens y)=qy\tens x,
\quad \Psi(y\tens x)=qx\tens y+(q^2-1)y\tens x}
extended to products by functoriality. This is the braiding induced
by the action of $\tilde{U_q(su_2)}$ on the quantum-braided plane
and is usually given by an array of $su_2$-type R-matrices,
see\cite{Ma:poi}\cite{Ma:book}.

\begin{propos} The restricted dual-quasitriangular comeasuring
bialgebra $M(R,\C_q^2)$
is generated by $1,s_{(i,j)},t_{(i,j)}$, $i,j\in\Z_+$ with the
relations in Proposition~3.1 and the additional relations
\[ R^{(i,j)}{}_{(a,b)}{}^{(k,l)}{}_{(c,d)}s_{(a,b)}s_{(c,d)}
=q^2 s_{(k,l)}s_{(i,j)},\quad
R^{(i,j)}{}_{(a,b)}{}^{(k,l)}{}_{(c,d)}t_{(a,b)}t_{(c,d)}=q^2
t_{(k,l)}t_{(i,j)}\]
\[ R^{(i,j)}{}_{(a,b)}{}^{(k,l)}{}_{(c,d)}s_{(a,b)}t_{(c,d)}
=q t_{(k,l)}s_{(i,j)},\quad
R^{(i,j)}{}_{(a,b)}{}^{(k,l)}{}_{(c,d)}t_{(a,b)}s_{(c,d)}=q
s_{(k,l)}t_{(i,j)}+(q^2-1)t_{(k,l)}s_{(i,j)}\] where $R$
corresponds to the braiding $\Psi$ on $\C_q^2$. Summation over
$a,b,c,d\in\Z_+$ should be understood.
\end{propos}
\proof We already know $M(\C_q^2)$ from Section~3.7, and now
quotient this further
by the $R\vect_1\vect_2=\vect_2\vect_1R$ relations. On the other
hand, the braiding on the generators (\ref{xybra}) immediately
gives
$R^{(a,b)}{}_{(1,0)}{}^{(c,d)}{}_{(1,0)}=\delta^a_1\delta^b_0
\delta^c_1\delta^d_0
q^2$  corresponding to $\Psi(x\tens x)=q^2 x\tens x$, etc. Thus the
additional relations on the generators $s_{(i,j)},t_{(i,j)}$ reduce
to the four sets as shown. They can be stated more compactly as the
relations of a rectangular\cite{MaMar:glu} quantum matrix
$A(R:R_{su_2})$, where $R_{su_2}$ is the standard $su_2$-type
$R$-matrix.
\endproof

The relations of $M(R,\C_q^2)$ may be further expressed in terms of
$R_{su_2}$ or, alternatively, defined by induction. If $x_i=(x,y)$
is an $SU_q(2)$-covariant covector notation then
\[ x_{i_1}x_{i_2}\cdots x_{i_n}=q^{p(i_1,\cdots,i_n)}x^{\#(i_1,
\cdots,i_n)}y^{n-\#(i_1,\cdots,i_n)}\]
defines $\Z_+$-valued functions $p,\#$ (here $\#$ is the number of
indices with value 1, and $p$ is the number of indices with value
$1$ to the right of each index with value $2$.) Then
\[ R^{(a,b)}{}_{(i,j)}{}^{(c,d)}{}_{(k,l)}=\delta^{a+b}_{i+j}
\delta^{c+d}_{k+l}
\sum_{\#(b_i)=c}\sum_{\#(a_i)=a}q^{p(b_i)+p(a_i)} Z^{a_1\cdots
a_{i+j}}{}_{1\cdots 2\cdots}{}^{b_1
\cdots b_{k+l}}{}_{1\cdots 2\cdots}\]
where there are $i$ indices $1\cdots$, $j$ indices $2\cdots$,
followed by $k$ indices $1\cdots$, and $l$ indices $2\cdots$, and
where $Z$ is the `partition function' array of copies of $R_{su_2}$
corresponding to $\Psi$ in the braided covector description of the
quantum-braided plane (see \cite[Thm. 10.2.1]{Ma:book}).
Alternatively, the general $\Psi$ in our basis $e_{(i,j)}=x^iy^j$
may be obtained by induction, via the formulae
\[ \Psi(y^i\tens x^j)=x\Psi(y^i\tens x^{j-1})q^i+(q^2-1)
[i]_{q^2}q^{2(j-1)} y\Psi(y^{i-1}\tens
x^{j-1})x\] and
\[ \Psi(x^iy^j\tens x^ky^l)=q^{i(2k+l-j)+l(2j-k)}y^l\Psi(y^j
\tens x^k)x^i.\]
These expressions follow from the functoriality (\ref{braalgpsi})
of the braiding.

\begin{corol} There is a bialgebra surjection $M_0(R,\C)\to M_q(2)$
and $M_q(2)$ appears as
a subalgebra (\ref{qmatst}) covered by this surjection. The
dual-quasitriangular structure of $M_0(R,\C)$ extends that of
$M_q(2)$.
\end{corol}
\proof Since $\Psi$ preserves the total degree, the only way to
obtain a non-zero coefficient of $x\tens x$, $x\tens y$, $y\tens x$
or $y\tens y$ in $\Psi(x^iy^j\tens x^ky^l)$ is with $i+j=1=k+l$.
Hence to compute $R^{(1,0)}{}_{(i,j)}{}^{(1,0)}{}_{(k,l)}$ etc., we
need only to consider (\ref{xybra}). Writing $1\equiv(1,0)$ and
$2=(0,1)$, the only nonzero entries of this form are given by the
standard $R_{su_2}$. Thus, the lowest level generators
$s_{(1,0)},s_{(0,1)}$ obey the relations of the quantum-plane in
R-matrix form, i.e. $ca=q ac$ in the notation (\ref{qmatst}).
Similarly, $t_{(1,0)},t_{(0,1)}$ form a quantum-braided plane, i.e.
$db=qbd$. These relations already hold in $M_0(\C_q^2)$. The third
relation in Proposition~4.3 similarly reduces to the four relations
$ba=qab$, $dc=qcd$, $cb=bc$ and $ad-da=(q^{-1}-q)bc$, while the
forth is then redundant. Hence the relations among these lowest
level generators of $M_0(R,\C_q^2)$ are precisely the relations of
the $2\times 2$ quantum matrices $M_q(2)$. \endproof

The geometric meaning of this is as follows. The surjection
corresponds classically to the inclusion of $2\times 2$ linear
transformations among the algebraic diffeomorphisms of the plane.
The inclusion of $M_q(2)$ corresponds classically to the projection
which associates to a diffeomorphism fixing zero its differential
at zero, i.e. the linear transformation induced on the tangent
space at zero.

Similarly, there is a surjection from $M(R,\C_q^2)$ to the
q-deformed Weyl algebra $\C_q^2\lbiprod M_q(2)$ cf \cite{Ma:book}
and an inclusion of it. One replaces the dilaton-extended
$SU_q(2)$, i.e. $GL_q(2)$, in the bosonisation construction
$\C_q^2\lbiprod GL_q(2)$ from \cite{Ma:poi}\cite{Ma:book} by the
quantum matrices $M_q(2)$. Conversely, at least these quotients of
$M(R,\C_q^2)$ become Hopf algebras by adjoining the inverse of the
$q$-determinant. Moreover, our constructions are quite general and
can be applied similarly to quantum-braided planes associated to
other $R$-matrices than $R_{su_2}$.

\subsection{Braided comeasuring bialgebras}

Finally, we give braided versions of all our comeasuring
constructions. From a categorical point of view we can fix any
braided category in which our algebra lives and look for the
universal comeasuring object in this braided category. The result
will now be a braided group\cite{Ma:exa} or bialgebra in the
braided category.

For the abstract aspect of these constructions, we use the
diagrammatic notation for `braided algebra' (due to the author) in
which we write products etc. as nodes
$\cdot=\epsfbox{prodfrag.eps}$ and `wire up' our algebraic
operation using $\Psi=\epsfbox{braid.eps}$,
$\Psi^{-1}=\epsfbox{braidinv.eps}$ as necessary. Note that the use
of `flow charts' to express algebraic or other operations is
nothing new -- it is used routinely by physicists as Feynman
diagrams and by engineers as wiring diagrams; the new feature in
braided mathematics \cite{Ma:bra}\cite{Ma:introp} is to adapt the
notation to algebraic constructions in braided categories, where
under and over crossings are nontrivial and distinct operators. We
use here the coherence theorem for braided
categories\cite{JoyStr:bra}.

In particular, using the braiding, one has a braided tensor product
algebra $A\und\tens B$ for $A,B$ algebras in the category. In a
concrete setting it is $(a\tens b)(c\tens d)=a\Psi(b\tens c)d$, but
more generally it is defined diagrammatically, see \cite{Ma:book}.
We define a comeasuring of $A$ as a pair $(B,\beta)$ where
$\beta:A\to A\und\tens B$ is a morphism and an algebra map.

\begin{propos} If $A$ is an algebra in a braided category, the
universal comeasuring object ${\und M}_1(A)$ is a braided group (a
bialgebra in the braided category).
\end{propos}
\proof This is shown in Figure~1. In part (a) we write the definition
of comeasuring in diagrammatic form. In part (b) we check that
comeasurings are closed under tensor product. The lower box is the
braided tensor product algebra $B\und\tens B$. Hence if $({\und
M},\beta_U)$ is the universal object, then $({\und M}\und\tens{\und
M},(\beta_U\tens\id)\beta_U)$ also comeasures, hence there is an
induced algebra map $\Delta:{\und M}\to{\und M}\und\tens{\und M}$.
Part(c) checks that it is coassociative.
\endproof
\begin{figure}
\[\epsfbox{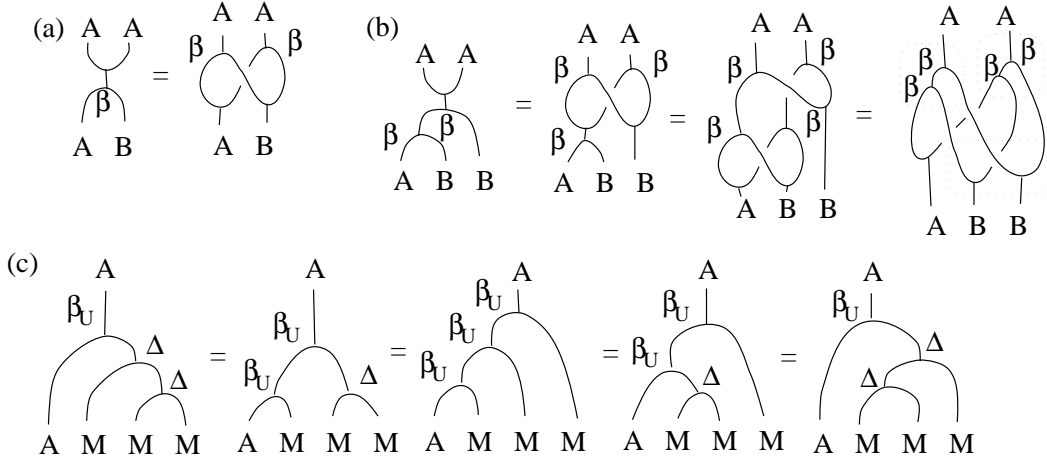}\]
\caption{(a) Comeasuring property of $(B,\beta)$ (b) Proof that
$B\und\tens B,(\beta\tens\id)\beta$
is another comeasuring (c) Proof that $\Delta$ is coassociative}
\end{figure}

As before, we have quotients ${\und M}(A)$ and ${\und M}_0(A)$ (the
latter when $1$ is split). We have explicit formulae in the
$R$-matrix case. For convenience we assume that $R$ is
bi-invertible in the sense that $\tilde{R}=((R^{t_2})^{-1})^{t_2}$
exists, where $t_2$ denotes transposition in the second matrix
factor.

\begin{propos} If $A$ is a finite-dimensional braided algebra with
braiding
determined by a biinvertible $R$-matrix, then the universal
comeasuring braided group has the explicit form ${\und M}_1(R,A)$
with generators $1$ and $u^i{}_j$ and the relations, coalgebra and
braiding
\[ c_{ij}{}^au^k{}_a=c_{ab}{}^k R^{-1}{}^a{}_c{}^b{}_d u^c{}_e
R^e{}_i{}^d{}_f u^f{}_j\]
\[ \Delta u^i{}_j=u^i{}_a\tens u^a{}_j,\quad \eps(u^i{}_j)
=\delta^i{}_j,
\quad \Psi(u^i{}_j\tens u^k{}_l)=u^m{}_n\tens u^r{}_s R^i{}_a{}^d{}_m
R^{-1}{}^a{}_r{}^n{}_b
R^s{}_c{}^b{}_l \tilde{R}{}^c{}_j{}^k{}_d.\] The braided coaction
and the braiding with $A$ is
\[ \beta_U(e_j)=e_a\tens u^a{}_j\quad \Psi(u^i{}_j\tens
e_k)=e_m\tens u^a{}_b R^{-1}{}^i{}_a{}^m{}_n R^b{}_j{}^n{}_k,\quad
\Psi(e_k\tens u^i{}_j)=u^a{}_b\tens e_m\widetilde R^n{}_k{}^i{}_a
R^m{}_n{}^b{}_j.\]
In compact form, the bialgebra structure is
\[ \vecu_3c_{12}{}^3=c_{12}{}^3R_{12}^{-1}\vecu_1R_{12}\vecu_2,
\quad \Delta\vecu=\vecu\tens\vecu,\quad\eps\vecu=\id,\quad
\Psi(R_{12}^{-1}\vecu_1\tens R_{12}\vecu_2)=
\vecu_2R_{12}^{-1}\tens\vecu_1R_{12}.\]
\end{propos}
\proof That we obtain here a bialgebra follows from the preceding
proposition
once we have established the universal property. Before doing this
we first outline, for completeness, the direct algebraic proof. For
this, we have to check that the coproduct extends as an algebra map
to the braided tensor product with the stated braiding $\Psi$.
Thus,
\align{&&\equad \Delta(c_{ab}{}^k R^{-1}{}^a{}_c{}^b{}_d u^c{}_e
R^e{}_i{}^d{}_f u^f{}_j)=
c_{ab}{}^k R^{-1}{}^a{}_c{}^b{}_d R^e{}_i{}^d{}_f  (u^c{}_p\tens
u^p{}_e) (u^f{}_q\tens u^q{}_j)\\ &&
=c_{ab}{}^k R^{-1}{}^a{}_c{}^b{}_d R^e{}_i{}^d{}_f
u^c{}_p u^m{}_n\tens u^r{}_s u^q{}_j R^p{}_u{}^z{}_m
R^{-1}{}^u{}_r{}^n{}_v R^s{}_w{}^v{}_q\tilde{R}{}^w{}_e{}^f{}_z\\
&&=c_{ab}{}^k R^{-1}{}^a{}_c{}^b{}_d  u^c{}_p R^p{}_u{}^d{}_m
u^m{}_n\tens u^r{}_s u^q{}_j  R^{-1}{}^u{}_r{}^n{}_v
R^s{}_i{}^v{}_q\\ &&=c_{dn}{}^au^k{}_a\tens u^r{}_s
R^{-1}{}^u{}_r{}^n{}_v u^q{}_j R^s{}_i{}^v{}_q
=c_{ij}{}^au^k{}_b\tens u^b{}_a= \Delta c_{ij}{}^a u^k{}_a.}
Or in the compact notation, this proof reads
\align{&&\equad \Delta( c_{12}{}^3R_{12}^{-1}\vecu_1 R_{12}\vecu_2)
=c_{12}{}^3R_{12}^{-1}
(\vecu_1\tens\vecu_1)R_{12}(\vecu_2\tens\vecu_2)=c_{12}{}^3
R_{12}^{-1}
\vecu_1\Psi(\vecu_1)R_{12}\tens\vecu_2)\vecu_2\\
&&=c_{12}{}^3R_{12}^{-1}
\vecu_1R_{12}\vecu_2\tens R_{12}^{-1}\vecu_1R_{12}\vecu_2=\vecu_3
c_{12}{}^3\tens R^{-1}_{12}\vecu_1
R_{12}\vecu_2=\vecu_3\tens\vecu_3 c_{12}{}^3=\Delta\vecu_3
c_{12}{}^3.}

We also have to  check that $\Psi$ is itself well-defined when
extended to products by functoriality (i.e. such that ${\und
M}_1(R,A)$ is a braided algebra). Here a direct proof is too
complex to write out in explicit terms and we give it only with the
compact notation.  Thus,
\align{&&\equad\Psi(\vecu_1R_{12}\tens c_{34}{}^2R_{34}{}^{-1}\vecu_3
R_{34}\vecu_4)=
\Psi(\vecu_1\tens  c_{34}{}^2R_{14}R_{13}R_{34}{}^{-1}\vecu_3
R_{34}\vecu_4)\\
&&= c_{34}{}^2R_{34}^{-1}\Psi(\vecu_1\tens
R_{13}R_{14}\vecu_3R_{34}\vecu_4)=c_{34}{}^2R_{34}^{-1}R_{13}
\vecu_3R_{13}^{-1}
\Psi(\vecu_1R_{13}\tens R_{14}R_{34}\vecu_4)\\
&&=c_{34}{}^2R_{34}^{-1}R_{13}\vecu_3R_{13}^{-1}
R_{34}\Psi(\vecu_1\tens
R_{14}\vecu_4)R_{13}=c_{34}{}^2R_{34}^{-1}R_{13}\vecu_3R_{13}^{-1}
R_{34}R_{14}\vecu_4
R_{14}^{-1}\tens\vecu_1 R_{14}R_{13}\\
&&=c_{34}{}^2R_{34}^{-1}R_{13}R_{14}\vecu_3R_{34}R_{13}^{-1}\vecu_4
R_{14}^{-1}\tens\vecu_1
R_{14}R_{13}=c_{34}{}^2R_{14}R_{13}R_{34}^{-1}\vecu_3R_{34}\vecu_4
R_{13}^{-1}
R_{14}^{-1}\tens\vecu_1 R_{14}R_{13}\\ &&=
R_{12}c_{34}{}^2R_{34}^{-1}\vecu_3R_{34}\vecu_4R_{13}^{-1}
R_{14}^{-1}\tens\vecu_1 R_{14}R_{13}=
R_{12}\vecu_2c_{34}{}^2R_{13}^{-1} R_{14}^{-1}\tens\vecu_1
R_{14}R_{13}\\ &&=R_{12}\vecu_2\tens
R_{12}^{-1}c_{34}{}^2\vecu_1R_{14}R_{13}=R_{12}\vecu_2\tens
R_{12}^{-1}\vecu_1 R_{12}c_{34}{}^2=\Psi(\vecu_1 R_{12}\tens
\vecu_2 c_{34}{}^2)} using repeatedly the QYBE and one of the two
covariance conditions (\ref{braalgmat}). On the other side, we have
\align{&&\equad\Psi(c_{12}{}^3R_{12}^{-1}\vecu_1R_{12}\vecu_2\tens
R_{14}R_{24}\vecu_4)
=c_{12}{}^3R_{12}^{-1}\Psi(\vecu_1R_{12}R_{14}\vecu_2\tens
R_{24}\vecu_4)\\
&&=c_{12}{}^3R_{12}^{-1}\Psi(\vecu_1 R_{12}R_{14}R_{24}\tens
\vecu_4)R_{24}^{-1}\vecu_2R_{24}=c_{12}{}^3R_{12}^{-1}R_{24}
\Psi(\vecu_1 R_{14}\tens
\vecu_4)R_{12}R_{24}^{-1}\vecu_2R_{24}\\
&&=c_{12}{}^3R_{12}^{-1}R_{24}R_{14}\vecu_4R_{14}^{-1}\tens\vecu_1
R_{14}R_{12}R_{24}^{-1}\vecu_2R_{24}
=c_{12}{}^3R_{14}R_{24}R_{12}^{-1}\vecu_4R_{14}^{-1}\tens R_{24}^{-1}
\vecu_1 R_{12}R_{14}\vecu_2R_{24}\\
&&=R_{34}\vecu_4 c_{12}{}^3 R_{12}^{-1}R_{14}^{-1}R_{24}^{-1}\tens
\vecu_1 R_{12}R_{14}\vecu_2R_{24}=R_{34}\vecu_4
c_{12}{}^3 R_{24}^{-1}R_{14}^{-1}R_{12}^{-1}\tens \vecu_1
R_{12}R_{14}\vecu_2R_{24}\\ &&=R_{34}\vecu_4
R_{34}^{-1}c_{12}{}^3R_{12}^{-1}\tens \vecu_1
R_{12}R_{14}\vecu_2R_{24}=R_{34}\vecu_4 R_{34}^{-1}\tens
\vecu_3c_{12}{}^3 R_{12}R_{14}\\
&&=R_{34}\vecu_4R_{34}^{-1}\tens\vecu_3R_{34}c_{12}{}^3=\Psi(\vecu_3
R_{34}\tens\vecu_4c_{12}{}^3)
=\Psi(\vecu_3c_{12}{}^3\tens
R_{14}R_{24}\vecu_4)} using the QYBE and the other half of
(\ref{braalgmat}). The proof for higher products is similar and the
general case follows by induction.

Finally, we verify the universal property. Thus, if $B$ is a
braided comeasuring on $A$, we define $\pi(u^i{}_j)\in B$ as
$\beta(e_j)=e_a\tens\pi(u^a{}_j)$. That this is a morphism implies
that the braiding with $A$ may be computed before or after applying
$\beta$. This is shown in Figure~2. Thus
\[ e_a\tens\Psi(\pi(u^a{}_i)\tens e_j)=\Psi^{-1}(e_b\tens e_c)\tens
u^c{}_a R^a{}_i{}^b{}_j
=e_m \tens e_n R^{-1}{}^m{}_c{}^n{}_b\tens u^c{}_a R^a{}_i{}^b{}_j\]
which tells us that the braiding $\Psi:B\tens A\to A\tens B$ is
compatible via $\pi$ with the braiding of ${\und M}(R,A)$ with $A$.
The braiding with $A$ on the other side similarly gives
compatibility of $\Psi:A\tens B\to B\tens A$.  Once these are known
then the braiding of $A$ with $\pi(u^j{}_k)\in B$ before and after
$\beta$ yields compatibility with $\Psi:B\tens B\to B\tens B$. This
is shown on the right in Figure~2. In this way, the braiding on
$\pi(u^i{}_j)$ has the same form as the braiding on $u^i{}_j$.
Next, the assumption that $\beta$ is an algebra map to $A\und\tens
B$ translates as
\[\beta(e_ie_j)=c_{ij}{}^a e_k\tens\pi(u^k{}_a)=(e_a\tens\pi(u^a{}_i))
(e_b\tens\pi(u^b{}_j))
=e_ae_c\tens \pi(u^d{}_e)R^{-1}{}^a{}_d{}^c{}_f R^e{}_i{}^f{}_b
\pi(u^b{}_j).\]
Thus $u^i{}_j\mapsto \pi(u^i{}_j)$ extends as an algebra map
$\pi:{\und M}(R,A)\to B$ and a morphism in the braided category.
\endproof
\begin{figure}
\[\epsfbox{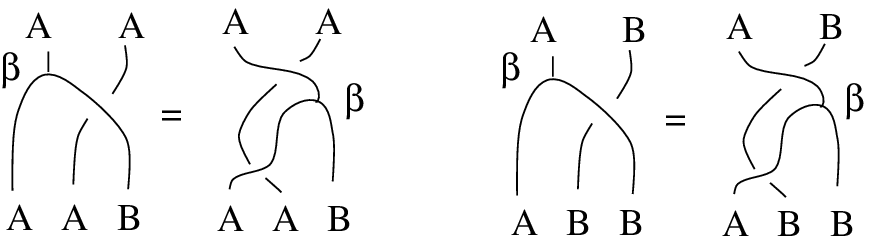}\]
\caption{Construction of braidings $B\tens A\to A\tens B$ and hence
$B\tens B\to
B\tens B$}
\end{figure}

Finally, we recall that the quantum matrices $A(R)$ above have a
braided group version $B(R)$\cite{Ma:exa} defined by a matrix of
generators $\vecu=(u^i{}_j)$ and relations $R_{21}\vecu_1
R\vecu_2=\vecu_2 R_{21}\vecu_1 R$, forming a braided group with the
above matrix coalgebra and braiding. The algebra relations here are
known in other contexts too\cite{FRT:lie} and sometimes called
`reflection equations', although they have been introduced and
studied as quadratic algebras by the author under the heading of
the braided matrix relations. Key properties such as its covariance
properties (as a braided algebra), the braided coproduct, results
about the Poincar\'e series etc. were introduced in \cite{Ma:exa}.

On the other hand, the quantum and braided matrices are closely
tied by a theory of transmutation which relates their products,
e.g.
\eqn{trans}{ \vect=\vecu,\quad R_{12}\vect_1\vect_2=\vecu_1R\vecu_2.}
See \cite{Ma:book}\cite{Ma:introp} for an introduction to this
transmutation theory of braided groups.

\begin{corol} The transmutation of $M_1(R,A)$ is the braided
comeasuring bialgebra
${\und M}_1(R,A)$ generated by $1$ and $u^i{}_j$ with the braided
matrix relations $R_{21}\vecu_1 R\vecu_2=\vecu_2 R_{21}\vecu_1 R$
and the further relations $\vecu_3c_{12}{}^3=c_{12}{}^3
R^{-1}_{12}\vecu_1R_{12}\vecu_2$ and coalgebra of ${\und M}_1(A)$
from the preceding proposition.
\end{corol}
\proof This is immediate from (\ref{trans}) applied to the relations
$\vect_3 c_{12}{}^3=c_{12}{}^3\vect_1\vect_2$ of $M(R,A)$.
\endproof

One may then proceed to apply the extensive theory of braided
groups to ${\und M}_1(A)$, ${\und M}_1(R,A)$ and their unital
quotients. For example, associated to any braided group in the
category of comodules of a quantum group $H$ (which is essentially
the situation above, with $H$ obtained from $A(R)$), one has
ordinary bialgebras ${\und M}_1(A)\lbiprod H$ etc., by the
bosonisation construction.

As a simple example, we consider $\C[x]$ as a braided algebra with
R-matrix (\ref{braline}). In this case the relations from
Proposition~4.6 are
\[ u^k{}_{i+j}=\sum_{a+b=k} u^a{}_i u^b{}_j q^{(i-a)b}.\]
Hence ${\und M}(\C[x])=\C<u_i|i\in\Z_+>$ is the free algebra as in
Section~3.1, but the other generators (and hence the matrix
coalgebra) are given from these by
\[ u^i{}_j=\sum_{a_1+\cdots +a_j=i} u_{a_1}\cdots
u_{a_j}q^{-(i-a_1)-\sum_{s=2}^i(a_1+\cdots a_{s-1})a_s}\], which is
a q-deformation of (\ref{deltapoly} as a braided group. The
braiding is
\[ \Psi(u_i,u_j)=u_j\tens u_i q^{(i-1)(j-1)}.\]
Similarly, the braided-commutativity in Corollary~4.7 in this
simplest example reduces to the usual commutativity relations. For
example, ${\und M}_0(R,\C[x])=\C[u_i|\ i\in\N]$ is the same algebra
as $\Diff_0(\C[x])$ in the classical case, but the coalgebra is
$q$-deformed and provides a nontrivial braided group structure on
this algebra.


\end{document}